\documentclass[12pt]{article}
\usepackage{amsmath}
\usepackage{amssymb}
\usepackage[utf8x]{inputenc}
\usepackage{enumerate}
\usepackage{makeidx}

\usepackage[pdftex,bookmarksnumbered,%
colorlinks,backref,pagebackref,hyperindex]{hyperref}
\hypersetup{pdftitle=THE SIGNATURE OF A MANIFOLD,
pdfauthor=Jose del Carmen Rodriguez Santamaria,
pdfsubject=differential geometry,
pdfkeywords=signature;cohomoly;index theorem}

\makeindex

\newtheorem{teo}{}

\begin{document}

\title{THE SIGNATURE OF A MANIFOLD}         
\author{Jos\'e Rodr\'\i guez}        
\date{}          
\maketitle

\begin{center}

\tableofcontents

\small{Abstract}

\end{center}

\normalsize

Let us consider   a compact oriented  riemannian manifold $M$ without boundary
and of dimension $n=4k$. The signature of $M$ is defined as the signature of a
given quadratic form $Q$. Two different products could be used 
to define $Q$ and
they render equivalent definitions: the exterior product of $2k$-forms and the cup
product of cohomology classes. The signature of a manifold is proved to yield a
topological invariant. Additionally, using the metric, a suitable Dirac operator
can be defined whose index coincides with the signature of the manifold. This
second version includes corrections  and  many examples.

\section{Introduction}

In differential geometry we can integrate  forms and only  forms.  The degree of
the form is given by the dimension of the domain of
integration. Say, if we have a three dimensional  density of electric charge,
the total charge is its integral: the charge density shall be encoded as a
3-form. Likewise, if we want to reproduce the fact that the line integral of a
force gives a work, we encode forces as 1-forms. On the other hand, Force =
Electric field x charge, so that the electric field could be encoded as a 1-form
also. The Lorentz force represents the capability of the electromagnetic field
to modify the world: it is a 1-form given by $\vec F=k\vec E+c\vec B\times \vec
v$, where $\vec E$ is the electric field, a 1-form, and  $\vec B$ is the
magnetic field that must be encoded as  a 2-form that operates over
a pair of
vectors, velocity $\vec v$ and  displacement.  

\

Let us imagine now a rubber surface    filled in electric charges with a given
density. Clearly, if we deform the surface, the total electric charge is
conserved. Thus, thanks to electric charges we have an invariant. The question
is whether or not we can somehow balance the role of the densities of charge in
such a way that we get a charge-free invariant, i.e., an invariant that depends
only on the manifold. In the following we present a positive partial answer to
this question at whose root we find cohomology, which is a byproduct of
derivation. 

\begin{teo}
\textbf{Differential forms}. The machinery for derivation and integration in differential geometry 
is built around forms that 
are defined over a smooth manifold $M$. 

Let $T_pM$ be the tangent space to $M$ at $p$. 
At any point $p \in M$, a $k$-form $\beta$ defines an alternating multilinear map from  $k$ factors 
of  $T_pM$ into $\mathbb{R}$:

    $$\beta_p\colon T_p M\times \cdots \times T_p M \to \mathbb{R}$$

An element in the domain can be understood as a parallelepiped localized at $p$ and the form 
measures its content, say, of electric charge.  

The set of all differential $k$-forms on a manifold $M$ at $p$ is a vector space  denoted 
$\Lambda^k_p(M)$.
The set of all differential $k$-forms on a manifold $M$ is a vector space  denoted $\Omega^k(M)$. 
So, a charge density is a 3-form, i.e.,  an 
element of $\Omega^3(M)$.   
The space of  all forms is denoted as $\Omega(M)$ (Nakahara,  \cite{Nakahara90} 
1990;  Frankel,  \cite{Frankel97} 1997).
\end{teo}

\begin{teo}
\textbf{Definition. } Let us consider a closed riemannian manifold with a metric $(M,g)$. 
The  \textbf{inner product of the two   $p-$forms},
$\omega, \eta$,   is defined as

\

$(\omega, \eta) = (1/p!)\int_M
\omega_{\mu_1...\mu_p}\eta^{\mu_1...\mu_p}\sqrt{g}dx^1...dx^n $.

\

\end{teo}

\begin{teo}
\textbf{Definition. } The \textbf{Hodge-star operator} \textbf{*}
 associates to
anyone $p-$form $\eta$
the one and only $(n-p)$-form $*\eta$  such that for whatever $p-$form $\omega$
we have:

\

$(\omega, \eta) =    \int \omega
\wedge *\eta $.

\

\end{teo}

\begin{teo}
\textbf{Derivatives}. Let $d_k$ denote the exterior derivative of $k-$forms and
$\delta_{k+1}$ its adjoint according to the scalar product $(\alpha , \beta ) =
\int_\Omega \alpha \wedge * \beta$.  We have that $d_{k+1} \circ  d_k = 0$. The
subindex of the exterior derivative rarely is explicitly written and $d_k$ is
simplified into $d$ and $d_{k+1} \circ  d_k = 0$ is  noted simply as $d^2=0$.
The same abuse of notation will be applied to other operators that appear in the
sequel. A form $\omega $ is called closed if $d\omega =0$, coclosed if $\delta
\omega =0$, exact if there exist $\psi$ such that $d\psi = \omega$, coexact if
there exists $\phi$ such that $\delta \phi = \omega$.
\end{teo}

We know that the divergence of a rotational is zero, and that the rotational of
a gradient is also zero. These facts lead to the general law $d^2=0$ which more
exactly means $d_{k+1} \circ  d_k = 0$.

\

\begin{teo}
\textbf{Cohomology. } In cohomology we take closed forms to investigate whether
or not they are exact.
\end{teo}

Since $d_{k+1} \circ  d_k = 0$ then $Im(d_k)\subseteq Ker(d_{k+1})$ and moreover
both are vector spaces.
The quotient $H^{k}=Ker(d_{k+1})/Im(d_k)$ is a vector space that is called the
\textbf{de Rham k-cohomology vector space}.

\

How are the classes of $H^k$? If  $\omega_1$ and $\omega_1'$ belong to the same
cohomology class then they
belong in the kernel of $d_{k+1}$, therefore $d_{k+1}\omega = d_{k+1}\omega'=0$
and moreover they differ by an exact form  $d_k\tau$, so  $\omega_1 = \omega_1'
+ d_k\tau$ with $d_{k+1} \circ  d_k \tau =0$. On the other hand, the operator
$d$ is antisymmetric.
Hence, $d_{k}\omega =0$ for $k>n$ given that $n$ is the dimension of the
manifold. This implies that $H^{k}(M)=0$ for $k>n$.   

\

We have presented cohomology in terms   of classes of closed
forms modulo exact differential forms, that operate locally. But cohomology can
also be given in terms
of holes, great topological items. For instance, to say that $H^{k}(M)=0$ for
$k>n$ is equivalent to saying that holes cannot have a dimension higher than
$n$.

\begin{teo}
\textbf{Example: the circle vs.  $\mathbb{R}$}.
\end{teo} 

Let us consider $\mathbb{R}$ and a
closed  smooth trajectory $\gamma $ over it that begins and ends in $a$ together
with a 1-form $\omega = f(x)dx$ with $f$ a continuous function. Then, there
exist $F$ such that $dF = f$ and $\int_\gamma w = \oint_a^a dF = F(a)-F(a) =
0$. By contrast, let us turn now to the circle. It has a problem: trajectories
appear that are closed and smooth but that does not un-walk what was walked. So,
the form that measures arc length,  $d\theta$, has the property : $\oint d\theta
\ne 0$. This seems to contradict the proposition saying that $\oint_a^a d\theta
= \theta(a) -\theta (a) = 0$. The veto to this and  similar reasonings  about
all
kinds of holes in higher dimensions is cohomology. This mathematical concept
is an abstraction that is embedded in our world: the existence
of electric
motors is due to the
fact that $H^1(\mathbb{R}^2 -\{0\})$ is not trivial. Similar effects are
observed in the theory of fundamental interactions. 

\

The cohomology of the circle in terms of differential forms reads  as follows: 
differential forms of degree 0 consist on functions $h(\theta)$ that take on
real values.
  A real
function $f$ is closed if $df = 0$, i.e., if $f$ is constant. Constant functions
= $H^0$ forms a vector space of dimension 1. At the other hand, Differential
forms of degree 1 are of the form   $\omega = f(\theta) d\theta$. They  are
always closed because $d\omega = (\partial
f/\partial \theta ) d\theta \wedge d\theta = 0$. If one has two differential
closed forms 
$\omega_1 = f(\theta) d\theta$ and $\omega_2 = g(\theta) d\theta$ then, 
 $\omega_2 = \omega_1  + (g(\theta) - f(\theta)) d\theta =  \omega_1  +
dF  $ if $f$ and $g$ are continuous. Therefore $H^1$ also has dimension 1 and
is generated by arc-length or by whatever other not null 1-form.

\

The general duality among differential forms and holes is consecrated by  the  Stokes' theorem.  
Thinking of the electric field would help us to understand the situation: if an electric field has a 
net flux across a closed surface that borders a spatial region $M$, it is because inside $M$ there 
shall be   electric charges that function like sources or sinks of field. So, the flux that is an 
integral over the border of $M$   must be equal to the net balance of  creation vs. annihilation  
of field inside $M$ that is related to the integration over $M$ of charge densities. Now, an 
electric  charge represents a hole in the domain of   the field  because 
there is no manner of defining it  to extend continuity and 
differentiability.  Formally:

\begin{teo} 
\textbf{Stokes' theorem. }  If a field can be represented by a  differential form $\omega$, then 
the net flux  of the field across the border $ \partial M $  of an orientable  manifold $M$  must 
be equal to the net    birth-death balance of field that happens at the interior of $M$ and that is 
measured by  $d\omega$:

 $$\int_{\partial M} \omega = \int_M d\omega  $$

The theorem literally states a truth for a differential form $\omega$ but
indeed it is also  true for its cohomology class:

\

$\int_{\partial M}[ \omega] =  \int_M d[\omega]$

\

\end{teo}

Proof:

\

 $\int_{\partial M}[ \omega] = \int_{\partial M} \omega + d\phi =  $
 $ \int_M
d(\omega  + d \phi) =  \int_M
d[\omega] =  \int_M
d\omega  + d^2 \phi =   \int_M
d\omega$. 

\

\begin{teo}
\textbf{Poincaré's lemma and Betti numbers. } 
\end{teo}

We use to understand cohomology in the light of the \textbf{Poincaré's lemma}:  
a closed form that is defined over a domain that is contractile into a point 
is also exact. Thus, a differential form operates locally but its exactness
depends  on large topological properties. In that way, cohomolohy connects local
and global properties. In fact, Poincaré's lemma allows us to see  the de Rham
cohomology vector space as an obstruction to the global exactness of closed
forms. The dimension of a
cohomology vector space is finite and is known as the \textbf{Betty number} of
the vector
space and
therefore could be seen as a measure of the variability of the global
inexactness of closed forms.

Thus, Betti numbers  measures obstructions to contractibility to a point, i.e., 
holes and disconnections in the domain:

\begin{enumerate}
 \item  $b_0$ is the number of connected components.
 \item    $b_1$ is the number of one-dimensional or \textquotedblleft circular"
holes.
 \item    $b_2$ is the number of two-dimensional holes or \textquotedblleft
voids".
\item $b_n$ is the number of n-dimensional holes.
\end{enumerate}

\

A sphere $S$ has one large n-hole and nothing else. So, all Betti numbers are
zero except $b_0 $ (which is 1 for $n > 0$ and 2 for $n=0$ because $S^0$
consists in two points) and $b_n = 1$.

\

Besides, we have the following three statements:

\begin{enumerate}
 \item if $M$ is contractile to a point then, by Poincaré's lemma, then all
closed forms are exact. Therefore, $H^p(M; \mathbb{R}) = 0$.
 \item If $M$ is a compact, connected, orientable manifold, and $dim M = n $,
then $H^n(M; \mathbb{R}) = \mathbb{R}$. This applies to all spheres: $H^n(S^n;
\mathbb{R}) = \mathbb{R}$.
\item If  $M$ is a compact, connected, non-orientable manifold, and $dim M = n
$, then $H^n(M; \mathbb{R}) = 0$. For $M = $  Möbius strip, $H^2(M; \mathbb{R})
= 0$.
\item If  $M$ is  a non-compact, connected, orientable or  non-orientable
manifold, and $dim M = n $, then $H^n(M; \mathbb{R}) = 0$. All $\mathbb{R}^n$
have  $H^n(\mathbb{R}^n; \mathbb{R}) = 0$.
\item Forms that measure lengths, areas volumes, ...,  are not exact on a
compact manifold. 
\end{enumerate}

\

We will illustrate the theory with  some specific calculations over  spheres,
toruses  and projective
spaces. So, let us review some material about them.

\

\begin{teo}
\textbf{Our notation for spheres is as follows: } 
\end{teo}

The $n-sphere$ of $\mathbb{R}^{n+1}$   is $S^{n} = \{ x \in \mathbb{R}^{n+1} |
d(x,0) = 1\}$. $S^n$ is a compact orientable manifold of dimension $n$.

The $n-disc$ of $\mathbb{R}^n$ is $D^n = \{ x \in \mathbb{R}^{n} | d(x,0) \le
1\}$. The $n-disc$ $D^n$ is a compact orientable manifold of dimension $n$.

The frontier of $D^n $ is the $n-1-sphere$: $\partial D^n = S^{n-1}$. The
frontier has one less dimension. 

The $n$- cell $e^n$ is a set that is homeomorphic to the open disc $D^n - 
\partial D^n$. By definition, $e^0$ is a point.

\begin{teo}
\textbf{The algebraic topology decomposition. } Our intuition suffers too much
in spaces of higher dimensions, so we recur to algebraic trickery. We can see
how it functions if we pay attention to a very simple case: $S^2$.
\end{teo}

If we punch $S²$ at a point, we get the sphere without a point: it is an open
set that can be covered by $\mathbb{R}^2$ as a single chart. So, it is also
equivalent to the open cell $ e^2$. We say that the punctured $S^2$
retracts to $e^2$. 

We can  specify how to reconstruct $S^2$: we  invert the retracting process of
the punctured $S²$ towards $e^2$. We use the jargon: glue the border of $e^2$ to
the point. Algebraic topologists say:

$S^2 = $ $e^o \cup e^2$.

The same procedure is valid for  every sphere.

\

One can transmit the same information using the Poincaré polynomial that results
from combining Betti number, $b_k$,  as coefficient and the  order of cohomology
as
power $k$ in the polynomial:

$p(t) = \Sigma  b_k t^k$

For the sphere $S^2$:

$p(t) = 1 + t^2$.

This says that  constant functions span $H^0$, that the order one cohomology of
the sphere  is zero, since every laze over the sphere can retract to a point. At
last,   $S^2$ has the volume element that is closed but not exact. So, the order
2
cohomology of the sphere is not zero. Now, since the sphere has only one 2-hole,
its
$H^2$ must have only one generator.

\

We will use Poincaré polynomials as  Ansatz to ease calculations. To see
how this works, let us consider the Torus $T = S^1 \times S^1$: if we punch
it, a normal donut, we find two non contractible  circles, one horizontal and
the other vertical, and a two dimensional cover. So, for $S^1 \times S^1$ we
have 

$p(t) = 1 + 2 t + t^2$.

that says that the Torus results from gluing a point to two circles to a 2-cell
that serves as cover: 

$ T = e^0 \cup 2e^1 \cup e^2$. 

\

Let us observe now that we
can find that polynomial as
the square of  $q(t) = 1 + t$, the polynomial of $S^1$:

\

$p(t) = (1+t)^2 = 1 + 2t + t^2$.

\

So, our guess is that the polynomial associated to the cartesian
product of two (compact, without boundary, oriented) manifolds is the product
of the polynomials associated to the manifolds. 

\

\begin{teo}
\textbf{Example. } Cohomology of a cylinder.
\end{teo}

An open cylinder is a surface that  has two
1-forms in cylindrical coordinates: $d\theta$ and $dz$. The first one is a
generator of $H^1$ while the second is not because  all pure-z closed curves are
contractile to a point.

\

\begin{teo}
\textbf{Quotient spaces}.
\end{teo}

We can infer how to work with quotient spaces if we look at the circle: it
is a
manifold in its own but it can be seen as a quotient space: the 
circle is the winding of  $\mathbb{R}$ that is furnished with the
equivalence relation $x \sim y$ iff $ y = x + 2k\pi$. Let us
notice
that winding preserves orientability. To see this, let us imagine that
$\mathbb{R}$ is
a trajectory that winds  over the circle. One chooses a tangent vector at a
given
point   over $\mathbb{R}$ and then one notices that the corresponding winding
always
observes the same direction: the circle is orientable and its orientation can be
given by an orientation over $\mathbb{R}$.

\

The element of volume in the circle is $d\theta$ which can be seen as the
differential form of coordinate $\theta$, the argument over the circle, else as
the
coordinate of $\mathbb{R}$ in the base space of the quotient space
$\mathbb{R}/\sim$.
We have that  $d\theta$ over $\mathbb{R}$ is exact but over the circle is not
because
there it measures arc length  which is $2\pi r$ over the whole circle over a
path
with start and final points coincident. There is no contradiction because
$\theta$ to
be a (uni-valued) chart function needs to be restricted to the open set
$(0,2\pi)$
and cannot be extended any further.

\

\begin{teo}
\textbf{The decomposition of the real projective plane. }
\end{teo}

Real projective space $\mathbb{RP}^n$ is defined to be the space of all lines
through the origin in $\mathbb{R}^{n+1}$. Formally: $\mathbb{RP}^n$ is the
quotient space of
$\mathbb{R}^{n+1} − \{0\}$ under the equivalence relation $v \sim \lambda v$ for
scalars $\lambda ≠ 0 $. 

\

The algebraic topology decomposition of the real
projective space is found by induction.

\

$\mathbb{RP} = \mathbb{RP}^1$ can be seen as follows:

\begin{enumerate}
 \item $\mathbb{RP} \subset \mathbb{R}^{2}$ is the set of lines that pass
through the origin. This space has dimension one: every such a line can be
represented by the two points of intersection of the line with the circle, the
1-sphere $S = S^1 \subset \mathbb{R}^2$, but
with antipodal points identified. One can imagine it as the superior hemicircle
with its terminal antipodal  points   identified. ( So, $\mathbb{RP}$ looks like
$S$.)

\item That superior hemicircle can be smashed into a disc  $D^1 \subset
\mathbb{R}^1$ to get that $\mathbb{RP}$ is the quotient
space of a disc $D^1$ with terminal points  $∂D^1$ identified. We get a circle:
$\mathbb{RP}^1 = S^1 = e^0 \cup e^1 $ that means that $S^1$ is just the open interval $(0,1)$ 
joined to a point. The topologies of $S^1$ and $\mathbb{RP}^1$ are the same but the winding number 
is multiplied
by two in $\mathbb{RP}^1$ that means that while you make a round trip in $S^1$ you gives two round 
trips in $\mathbb{RP}^1$.

\end{enumerate}

If we pass to the next dimension, we get:

\begin{enumerate}
 \item $\mathbb{RP}^2 \subset \mathbb{R}^{3}$ is just the sphere $S^2 \subset
\mathbb{R}^3$
with antipodal points identified. One can imagine it as the superior hemisphere
with the antipodal  points of the bottom circle identified. (This circle
functions as a single point, so $\mathbb{RP}^2$ looks like $S^2$.)

\item That superior hemisphere can be smashed into a disc  $D^2 \subset
\mathbb{R}^2$ to get that $\mathbb{RP}^2$ is the quotient
space of a hemidisc $D^2$ with antipodal points of its border $S^1 = ∂D^2$
identified.
\item Since $∂D^2 = S^1$ with antipodal points identified is just
$\mathbb{RP}^{1}$, we see that $\mathbb{RP}^2$ is obtained  by attaching a
$2$-cell $e^2$ to $S^1 = \mathbb{RP}^{1}$. So, we get $\mathbb{RP}^2 = S^1 \cup e^2 = e^0 \cup 
e^1\cup e^2 $
\end{enumerate}

In general, we have:

\begin{enumerate}
 \item $\mathbb{RP}^n \subset \mathbb{R}^{n+1}$ is just the sphere $S^n$
with antipodal points identified. One can imagine it as the superior hemisphere
with the antipodal  points of the bottom hyper-circle identified.
\item This hemisphere can be smashed into a disc to get that $\mathbb{RP}^n$ is
the quotient
space of a hemidisc $D^n \subset \mathbb{R}^n$ with antipodal points of its
border $∂D^n$ identified.
\item Since $∂D^n$ with antipodal points identified is just $\mathbb{RP}^{n}$,
we see that $\mathbb{RP}^n$ is obtained  by attaching an
$n$-cell $e^n$ to  $\mathbb{RP}^{n-1}$ but taking care of gluing the structure
of lower dimension to the
border of the higher one.
\end{enumerate}

Therefore, the decomposition of  $\mathbb{RP}^n$ is  $e^0 \cup e^1 \cup ... \cup
e^n$
with one cell $e^i$ in each dimension $i ≤ n$. Because we have one hole in each
dimension, we might predict that all Betti number are one.  Nevertheless, one must take care of orientability:

\begin{enumerate}
 \item $\mathbb{RP} = S$ is orientable.
 \item $\mathbb{RP}^2 = S^2$ with antipodal points identified is not orientable:
if one walks over the sphere $S^2$ through a maximal circle and one carries in
parallel transport a frame $\vec i, \vec j$, then while one remains in the upper
hemisphere, the corresponding frame in $\mathbb{RP}^2$ has the same orientation.
But if one passes to the lower hemisphere, the corresponding frame in
$\mathbb{RP}^2$  reverses orientation and so we cannot have a single valued
orientation from the different parts of $S^2$: $\mathbb{RP}^2$ is not
orientable. Hence, its 2-Betti number is zero.

\item In general, $\mathbb{RP}^n$ is orientable for odd dimensions and
non-orientable for even dimensions.
\item The Betti numbers cannot be read from
the cell structure in even dimensions because this depends on gluing maps whose
output depends on orientability. 
\end{enumerate}

\begin{teo}
\textbf{Remark.}
$\mathbb{RP}^n$ can de understood as the compactification of
$\mathbb{R}^{n}$. To fix ideas, let us show that $\mathbb{RP}^2$ can be
seen as the compactification of the plane $\mathbb{R}^{2}$. To begin with, 
$\mathbb{RP}^2$ contains a copy of the plane $\mathbb{R}^{2}$.   In fact, every
point of the
plane $z= 1$ defines a unique line that passes through the origin of
$\mathbb{R}^{3}$: that is why we say that the plane $\mathbb{R}^{2}$ is
contained in $\mathbb{RP}^2$. Nevertheless, $\mathbb{RP}^2$ is bigger than the
plane $\mathbb{R}^2$. 
To get  $\mathbb{RP}^2$,  we need to add the set of horizontal lines. That set
is homeomorphic to  $S^1$ with antipodal points identified. Now, the
plane $z= 1$ is homeomorphic to the upper
part of $S^2$ which in its turn is homeomorphic to $e^2$, the open disc of the
plane $z=0$. So, the compactification of the the
plane $e^2$ is made by $S^1$ with antipodal points identified, which functions
as a point. At last, $\mathbb{RP}^2$ looks like a sphere.
\end{teo}

\begin{teo}
\textbf{The decomposition of the complex projective plane}.
\end{teo}

Complex projective space $\mathbb{CP}^n$ is the space of complex lines through
the origin in $\mathbb{C}^{n+1}$. Formally,   $\mathbb{CP}^n$ is   the quotient
space of $\mathbb{C}^{n+1} − \{0\}$ under the
equivalence relation $v \sim \lambda v$ for $\lambda \ne 0$. To fix ideas, let us consider the point $ P = (1,0,.. .,0) \in \mathbb{CP}^n$. The complex line through $P$ is the set composed of all points of the form $\lambda P = (\lambda,0,.. .,0) $ for $\lambda \in \mathbb{C}$. This set is just a copy of $\mathbb{C}$ of complex dimension 1 and real dimension 2. The same happens for every other line.

\

What is the essence of $\mathbb{CP}^n$ from a topological point of view? The
following construction shows us that $\mathbb{CP}^n$  can best be understood as
a generalized Riemann sphere with a recursive construction, i.e., as the
compactification of $\mathbb{C}^n$ that uses a simple constructive algorithm.
This results from  the way as one builds $\mathbb{CP}^n$ from
$\mathbb{CP}^{n-1}$:

\

\begin{enumerate}
\item    $\mathbb{CP}^n \subset \mathbb{C}^{n+1}$  has complex dimension $n$. A
line in $\mathbb{C}^n$ can be denoted by any non zero  point in it. Let us
denote the line in $\mathbb{C}^n$  through point $(z_1, ..., z_{n+1}) $ as $[(z_1, ...,
z_{n+1})] $. Let $j$ be one of those coordinates that are not zero. We have
$[(z_1, ... , z_{n+1})]  = [(\frac{z_1}{z_j}, ...,1,..., \frac{z_{n+1}}{z_j})]
\leftrightarrow [(w_1, ... , w_{n})]$ with $w_i \in \mathbb{C}^n$. We conclude
that $\mathbb{CP}^n$ can be covered with $n$ patches each one a copy of
$\mathbb{C}^n.$ We refer to this system of coordinates as homogeneous.

\item $\mathbb{CP}^n$ is compact. In aforementioned  coordinates $\mathbb{CP}^n$
looks  unbounded, so it is convenient to think that it inherits from  a sphere 
a  metric and a topology that makes it a compact manifold. In fact, the unit
sphere $S^{2n+1}$ of $\mathbb{C}^{n+1} = \mathbb{R}^{2n+2}$ is another cover of
$\mathbb{CP}^n$: every line that passes through the origin has a point of norm
one.  The corresponding
association  defines a continuous function from the sphere, which is compact,
onto $\mathbb{CP}^n$.
Since the image of a compact subset is compact, we have that $\mathbb{CP}^n$ is
compact. The topology
of $\mathbb{CP}^n$ is that of the quotient space  $S^{2n+1}$ by the
$\lambda$ relation. 

\item  $\mathbb{CP}^n$ can be understood as the compactification of 
$\mathbb{C}^n$. In first place, $\mathbb{C}^n$ can be embedded into
$\mathbb{CP}^n$ and, in second place, the points of $\mathbb{CP}^n$ that are not
in the embedding can be understood as the points in infinity of $\mathbb{C}^n$.
This is shown as follows. The embedding of  $\mathbb{C}^n$  into $\mathbb{CP}^n$
 is given by

 $(z_1, ..., z_{n}) \leftrightarrow [(1, z_1, ..., z_{n})] $

Equivalently, we observe $\mathbb{CP}^n$ trough the hyper-plane $z_o = 1 + 0i$ of real dimension $2n$.

Let us notice now that the points of the form $[(1, z_1, ..., z_{n+1})] $ do not
cover $\mathbb{CP}^n$ completely: we lack those points of the form $[(0, z_1,
..., z_{n+1})]$. Let us see now why these points can be understood as points at
infinite in  $\mathbb{C}^n$. In fact, we can represent a point of $\mathbb{C}^n$
in infinite  as

 $(\lambda z_1, ..., \lambda z_{n}) $ for $\lambda \rightarrow \infty$.

The embedding gives

$ [(1, \lambda z_1, ..., \lambda z_{n})]   =  [(\frac{1}{\lambda }, z_1, ..., 
z_{n})]$
$\rightarrow  [0, z_1, ...,  z_{n})]$

This shows that $\mathbb{CP}^n$ becomes compact by adding the points at
infinity of $\mathbb{C}^n$. Nevertheless, our procedure
seems to produce two discrete units. To remedy this trouble, we use the sphere
as intermediary: every  line that passes through the origin also can be
represented over the sphere $S^{2n+1}$. Actually, the demi-sphere on the side of
$(1,0...,0)$ suffices, whose border represents the points at infinity.
Concretely, the border is $S^{2n-1}$ because it is the sphere in $\mathbb{C}^n =
\mathbb{R}^{2n}$. So, this border with the $\lambda$ equivalence is what serves
for the compactification of 
$\mathbb{C}^n$.

\item The topological structure of $\mathbb{CP}^n$ is given by  $\mathbb{CP}^n =
e^0
\cup e^2 \cup  ... \cup e^{2n}$. To see this, let us observe that the points at
infinity $ [(0, z_1, ...,  z_{n})]$ are stable under elongation, as it should
be. But to be stable under elongation is the trade mark of projective spaces.
That is why the points at infinity represent an embedding of $\mathbb{CP}^{n-1}$ into $\mathbb{CP}^n$.
From
this we deduce the way to form   $\mathbb{CP}^n$ from $\mathbb{CP}^{n-1}$: we
must
paste $\mathbb{CP}^{n-1}$ to  the infinite of (the embedding of) 
$\mathbb{C}^n$,
which must be understood as an open set, which
in
its turn is homeomorphic to $e^{2n}$. By recursion we get $\mathbb{CP}^n = e^0 \cup
e^2
\cup  ... \cup e^{2n}$. This construction proposes that its   Poincaré polynomial
is

$p(t) = 1 + t^2 + ... + t^{2n}$.  

\

\end{enumerate}

\begin{teo}
\textbf{Remark.} Our representation of projective spaces by means of
quotient spaces of spheres  allows us to see that these spaces are bounded and
that every path contained in any one of them  always leads to   an interior 
point. So, projective spaces are closed and without boundary. Officially, they
are boundariless, compact manifolds. This representation also allows us to
associate a volume form and a volume to the entire manifold if it is orientable.
The volume of the manifold would be equal to just the volume  of the cell of
highest dimension contained in the given quotient space. Glued parts of lower
dimensions have measure zero. Since apart from our representation there exist
many others (just change the radius of the representing sphere), we will say
that the volume of the quotient manifold is an unspecified  positive number $c$
and our propositions will state relations with that number $c$.
\end{teo}

\begin{teo}
\textbf{Electric motors and cohomology.} 
\end{teo}

We know from experience that electromagnetism exercises forces over charged
particles that eventually can be in movement. A classical description of this
activity is given by the Lorentz force   $\vec F=k\vec
E+c \vec v  \times  \vec B$, where $E$ es the electric field,  $\vec B$ is the
magnetic field and  $\vec v $ is the velocity of the charged particle.  This
equation allows us to understand how an electric motor functions:

\begin{enumerate}
 \item The term $\vec F= c \vec v  \times  \vec B$ says that if the velocity of
an
electron is perpendicular to a magnetic field, then it is deflected
sidewards. If instead of one electron we consider an electric current along a
wire, the deflection of electrons will be passed   to the wire because electrons
can move along the wire but as they try to leave the wire, a separation of
charges of different polarity occurs and so the wire is attracted and dragged by
the electrons. 
\item If the wire forms a square, it will begin to spin. In fact, if one
of
its sides deflects in one direction, the opposite side will deflect in contrary
direction because  corresponding currents have opposite directions.
\item The spinning force   can be made to endure forever if an appropriate
switching mechanism is endowed to cause the direction of current to reverse in
agreement with the position of the wire (Nave, 2013).
\end{enumerate}

Now, classic electromagnetism and electric motors cannot be considered without
cohomology. Let us see why.

There are many versions of electric motors but
all rely on magnets. Some depend on the magnetic field created by currents
along a wire, other depend on permanent magnets.

In regard with electromagnets,  magnetic field
lines form circles around a (straight, infinite) wire carrying an electric
current. The magnitude
of this field grows as one approaches the wire. That is why  the magnetic
field
cannot be defined directly over the wire because it would be  not univalued.
Therefore, the domain of definition of the magnetic field has a hole composed
by the wire which is equivalent to a hole in a plane.

 On the other side, permanent magnet results
from atoms that function as tiny magnets whose fields are aligned. One can
imagine in classical mechanics that these atoms have an unbalanced electron that
spins around the nucleus and that this movement generates a magnetic field.
That field is well defined everywhere with exception of the trajectory of the
electron: we end with the same cohomology: that of a plane with a hole.

\

The
next reasoning shows that the non nil cohomology  is not just a classical effect
but
that it is endemic to electromagnetism.

 In the depiction of electromagnetism given by the Lorentz force, 
it is a collage of two items: electricity and magnetism. The
unification of these two fields was given by relativity: electricity and
magnetism are two particular folds  of a single entity that populates
space-time,
the electromagnetic field.

The electromagnetic  field is described   by a 2-form $T$. It happens that $T$
can
be expressed as a differential: $T = dA$ where $A$ is a 1-form, the vector
potential. But, beware,  $T$ is an experimental
not nil object. This implies that $A$ itself   cannot be exact, i.e. there is
no scalar function $\alpha$ such that $A = d\alpha$. Otherwise $T = dA
= d^2\alpha = 0$. In other words, we are declaring that the cohomology vector
space $H^1$ is not nil. Therefore, the tensor field cannot have global
definition: it must have holes else be defined by sectors. De facto, the
electromagnetic field cannot be defined in the points that are occupied by 
electric charges. A charge defines a punctual hole in space but a line in
space-time. So, its cohomology in space-time is just that of a plane with a
punctual hole. That is why we say that the
existence of electric motors is due to the
fact that $H^1(\mathbb{R}^2 -\{0\})$ is not trivial.  

A
very soft introduction to these themes can be found in
Rodriguez ( \cite{Rodriguez08}, 2008).

\section{The signature quadratic form}

We define a bilinear form and find the constraints under which it becomes symmetrical.

\bigskip

\begin{teo}
\textbf{Example to show the whole idea. } 
\end{teo}

The Torus  $S^1
\times S^1$ is compact, without boundary,
orientable manifold of real dimension $2$. 
  Its  Poincaré polynomial  is:

$p(t) = (1+t)^2 =  1 + 2t + t^{2}$

This polynomial says that  a torus can be decomposed as a point to which two
circles are glued to
which a covering 2-cell must be attached. This means that  $H^1$ has 2 
generators, $d\theta$ and
$d\phi$,  while $H^2$ has 1, $d\theta d\phi$.

\

We can define a   
quadratic form as follows:

$$Q: H^{1}(M)\otimes H^{1}(M) \rightarrow \mathbb{R}$$

$$Q([\omega_1], [\omega_2])= \int_M \omega_1\wedge \omega_2 = \int
\omega_1\wedge \omega_2$$

\

This bilinear form can  be represented by a $2 \times 2$ matrix that we also call $Q$, whose entries 
represent respectively 

$\int d\theta \wedge d\theta = 0$, 

$\int d\theta \wedge d\phi = c$, in the same class as the area-form.

$\int d\phi \wedge d\theta = -c$,

$\int d\phi \wedge d\phi = 0 $

\

So,

$$Q =
\bordermatrix{~ &d\theta  & d \phi  \cr
                  d\theta & 0 & c \cr
                  d \phi  & -c & 0\cr}$$

\

The characteristic polynomial of this matrix  is:

$p(\lambda) = \lambda^2 + c^2. $

The roots of this polynomial are imaginary and cannot be compared with zero. The problem is just 
that matrix $Q$ is not symmetric, a failure that in its turn hangs  on the lack of 
commutativity of the wedge product. So, $Q$ will be a bilinear symmetric form when the wedge 
product is commutative, a property that holds in dimensions that    are multiple
of $4$.

\

\begin{teo}
\textbf{Definition. } Let $M$ be a compact oriented  riemannian manifold without boundary and of 
dimension n=4k. \textbf{The signature quadratic form} of $M$ is defined to be the  bilinear form

$$Q: H^{2k}(M)\otimes H^{2k}(M) \rightarrow \mathbb{R}$$

$$Q([\omega_1], [\omega_2])= \int_M \omega_1\wedge \omega_2 = \int
\omega_1\wedge \omega_2$$
\end{teo}

\begin{teo}
\textbf{Theorem.  } The signature quadratic form is well defined, i.e., the integral on the right
does not depend on the representatives of the cohomology classes.
\end{teo}

To prove claimed result,
let us take another representative of $[\omega_1]$ say $\omega_1' = \omega_1 +
d\tau$. Now, whatever representative one uses, one gets the same result:

\

$\int \omega_1'\wedge \omega_2= \int (\omega_1 + d\tau) \wedge \omega_2$
$=\int \omega_1  \wedge \omega_2 +\int  d\tau \wedge \omega_2$

$ \hspace{1.7cm}= \int \omega_1  \wedge \omega_2 +\int  d(\tau \wedge \omega_2)
$ 
$= \int \omega_1  \wedge \omega_2 +0 $

$ \hspace{1.7cm}= \int \omega_1  \wedge
\omega_2$

\

We have used the fact that $ d(\tau \wedge \omega_2) =  d\tau \wedge \omega_2
\pm \tau \wedge d\omega_2 = d\tau \wedge \omega_2$ because $\omega_2$ is closed.
And that 

\

$\int   d(\tau \wedge
\omega_2) = \int_{ M}  d(\tau \wedge
\omega_2) = \int_{\partial M}  \tau \wedge \omega_2 =0 $

\

because our manifolds have no boundary, like the spheres. Extending the reasoning to the general
form of representatives of $[\omega_2]$ we obtain the independence of the
bilinear form from representatives and so $Q$ is defined over pairs of cohomology classes
of order $2k$.

\

  Let us prove now that when $n = 4k$, we get a symmetric bilinear form:

\begin{teo}
\textbf{Theorem. } The signature quadratic form $Q$ is symmetric iff $n=4k$. Therefore, its 
eigenvalues are real.
\end{teo}
Proof: For  a $p$-form $\omega $   and  a $q$-form $\nu$   we have the general rule:

$\omega \wedge \nu = (-1)^{pq}\nu\wedge\omega$

For our case we have $p=q$ then $(-1)^{pq}= (-1)^{p^2}$. But $p^2=p$ mod 2 , i.e., $p^2-p= p(p-1) = 
2m$, i.e., 
either $p$ or $p-1$ is even and hence $p(p-1)$ is always even. So,  $(-1)^{p^2}=(-1)^{p}=1$ iff $p$ 
is even iff $2p$, the dimension of the manifold, is multiple of 4.

\

\begin{teo}
\textbf{Definition. } When $Q$ is symmetric, its eigen-values are real and can be compared with
zero. In that case, we can compute the \textbf{signature of the quadratic form}
$Q$  as the number of positive eigenvalues minus the number of negative ones. 
\end{teo}

\subsection{The cup product}

Since the wedge product is defined for pairs of forms of any order we can verify that  the vector 
space  $H^*= H^*(M)= \oplus_{p\ge 0}^n H^pM$ can be endowed with a ring structure as follows: for  
$[\omega] \in H^p$ and $[\nu] \in H^q$ we define the cup product as

$$\sqcup: H^* \times H^* \rightarrow H^*$$

$$[\omega] \sqcup [\nu]= [\omega \wedge \nu]$$

Let us check that this is a well defined product in the set of cohomology classes: if $[\omega] \in 
H^p$ and $[\nu] \in H^q$ then $\omega \wedge \nu \in \Omega^{p+q}$ and $[\omega \wedge \nu] \in 
H^{p+q}$ but we shall show that if we take other representatives their wedge product is still in 
$[\omega \wedge \nu]$.

If $\omega, \omega' \in [\omega]$ in $H^p$ then $d\omega= d\omega'=0$ and $\omega'= \omega+d\tau$. 
If $\nu \in H^q$ then $d\nu=0$. Therefore $\omega' \wedge \nu = (\omega +d\tau) \wedge \nu = \omega  
\wedge \nu  + d\tau \wedge \nu = \omega  \wedge \nu  + d(\tau \wedge \nu)$ because  $d(\tau \wedge 
\nu)=  d\tau \wedge \nu +  \tau \wedge d\nu = d\tau \wedge \nu$. Henceforth, $\omega'  \wedge \nu$ 
and $\omega  \wedge \nu$ differ by an element of the form $d(\alpha)$ where $\alpha = \tau \wedge 
\nu $ and so they belong to the same cohomology class.

\bigskip

Restricting the cup product to $H^{2k}(M)\otimes H^{2k}(M)$, the signature quadratic form $Q$ takes 
the form:

$$Q([\omega_1], [\omega_2])= \int [\omega_1]\sqcup [\omega_2] = \int \omega_1\wedge \omega_2$$

We see that that the cup product is
commutative when restricted to $H^{2k}(M)\otimes H^{2k}(M)$ because $Q$ is
symmetric when $n= 4k$.

\subsection{The signature of a manifold,  $\sigma(M)$}

 Since the space $H^{2k}$ is finite dimensional with dimension $b_{2k}$, called
the 2k-th Betti number, we can take a basis  $E= \{[h_{i}]\}$  with respect to
which  the quadratic form $Q$ has  an associated matrix, $Q_{[E]}$, whose
entries are $Q_{[E],ij}= \int [h_{*i}]\wedge [h_{*j}]$. This matrix can be
diagonalized to a matrix with eigenvalues $\lambda_1,
\lambda_2,...\lambda_{s}$., where $s= dimH^{2k}=b_{2k}$. Below we will see a
concrete basis in which Q is diagonal. 

\

Let us recall that the signature of
a real symmetric matrix is the number of positive eigenvalues minus the number
of negative eigenvalues. The signature of a matrix is invariant under  changes
of
basis with the same orientation. This happens because the signature of a real
symmetric
matrix measure an intrinsic property. In fact, a matrix represents a
linear transformation whose number of positive eigenvalues  is  the dimension of
the maximal vector subspace over which it is positive-definite.  All these
results enable the next

\bigskip

\begin{teo}
\textbf{Definition. }  The \textbf{signature of  compact, without boundary,
orientable manifold, $\sigma(M)$},  is the signature of the  quadratic form $Q$.
\end{teo}

Let us find the   signature  over some examples. To this aim, we will follow the general 
procedure stated above to calculate the signature of a riemannian boundariless
manifold and of dimension $4k$: it  is the number of positive eigenvalues minus
the number of negative eigenvalues of the finite dimensional matrix of $Q$ in
any basis of $H^{2k}$ and in any system of coordinates, where $Q$  is  the 
bilinear form

$$Q: H^{2k}(M)\otimes H^{2k}(M) \rightarrow \mathbb{R}$$

$$Q([\omega_1], [\omega_2])= \int \omega_1\wedge \omega_2$$

\begin{teo}
\textbf{Example. } Let us calculate the signature of $S^4$.
\end{teo}

$S^4$ contains only one big hole, so its decomposition    is

$S^4  = e^0   \cup e^{4} $

and its Poincaré polynomial is

$p(t) = 1 + t^4$.

There is no terms with intermediate powers because the sphere has no low
dimensional holes. 

Thus, to calculate the signature of $S^4$, we must consider the space $  H^{2}(S^4)\otimes 
H^{2}(S^4)  $ but  $H^2$  is
zero. So, the signature is zero.

\

\begin{teo}
\textbf{Example. } Let us find the signature  of  the
4-torus 

$ T^4 = S^1 \times S^1 \times S^1 \times S^1  $.
\end{teo}

We consider  $T^4$ as a compact  manifold over the reals of dimension 4.   Its 
Poincaré polynomial  is:

\

$p(t) = (1 + t)^4 = 1 + 4t + 6t^2 + 4t^3 + t^4 $   

\

or

\

 $T^4 = e^0 \cup 4 e^1 \cup 6 e^2 \cup 4e^3 \cup e^4$

 \

This  means that  $H^2$ has 6  generators while $H^4$ has 1. Thus, $Q$ is a $6
\times 6$ matrix. Specifically, $H^2$ is generated by $d\theta_1d\theta_2 $,
$d\theta_1d\theta_3 $, $d\theta_1d\theta_4 $, $d\theta_2d\theta_3 $,
$d\theta_2d\theta_4 $, $d\theta_3d\theta_4 $. On the other hand,
$d\theta_1d\theta_2 \theta_3\theta_4$is the generator of $H^4$ which is
in the same class as the form that measures
the 4-area of $T^4$. Let its integral  be $c$.   So the matrix $Q$ of the
integrals of wedge products is:

$$Q = \bordermatrix{~ & d\theta_1d\theta_2 & d\theta_1d\theta_3 &
d\theta_1d\theta_4 & d\theta_2d\theta_3 & d\theta_2d\theta_4 &
d\theta_3d\theta_4 \cr
d\theta_1d\theta_2 &                      0 & 0 & 0 & 0 & 0 & c \cr
 d\theta_1d\theta_3 &                     0 & 0 & 0 & 0 & -c & 0 \cr
 d\theta_1d\theta_4 &                     0 & 0 & 0 & c & 0 & 0 \cr
 d\theta_2d\theta_3 &                     0 & 0 & c & 0 & 0 & 0 \cr
  d\theta_2d\theta_4 &                    0 & -c & 0 & 0 & 0 & 0 \cr
  d\theta_3d\theta_4 &                    c & 0 & 0 & 0 & 0 & 0 \cr }$$

 For instance, $\int  d\theta_3d\theta_4 d\theta_1 d\theta_2 = $
 $-\int    d\theta_3d\theta_1d\theta_4      d\theta_2 =$
 $+\int   d\theta_1 d\theta_3d\theta_4     d\theta_2 =$
 $-\int    d\theta_1d\theta_3   d\theta_2   d\theta_4=$
 $+\int   d\theta_1d\theta_2     d\theta_3d\theta_4  = c $      
                  
\

The characteristic polynomial  is:

$p(\lambda) =  \lambda ^6 - 3 c^2 \lambda ^4 + 3 c^4 \lambda ^2 + c^6 = (\lambda - c) ^3  (\lambda 
+ c) ^3. $

This polynomial has six roots, three of them are  $c$ and the other 3 are  $-c$.

The number of positive eigenvalues equals that of negative ones.  Therefore, the
signature of   $T^4$ is zero.

\begin{teo}
\textbf{Example of a sophism about the  signature  of $\mathbb{CP}^2$.} The
clarification of this sophism will be found below.
\end{teo}

We consider $\mathbb{CP}^2$ as a compact  manifold over the reals of dimension
4. It   can be decomposed as $e⁰ + e^2 + e^4$. Its Poincaré polynomial  is:

$p(t) = 1 + t^2 + t^{4}$

which means that $b_2 = b_4 = 1$. So, $H^2$ and $H^4$ both have one generator.
Thus, $Q$ is a $1 \times 1$ matrix. Nevertheless, the only entry of $Q$ is
zero because every $Q$ matrix has zeros in its diagonal due to the fact that in
the diagonal appear   squares of forms that necessarily have repeated terms.
Therefore, $Q$ is the  zero matrix, a result that predicts that the only
eigenvalue of $Q$ is zero. Hence, the
signature
of   $\mathbb{CP}^2$ is zero.

\section{Invariance}

A diffeomorphism represents   a change of coordinates and we will show that the
signature of a manifold is independent of
them so, it is an intrinsic object. A direct demonstration of this result
over the ring $H^*$ is instructive and can   be  carried out thanks to the
pull-back technology. When in calculus one says  \textquotedblleft change of 
variable\textquotedblright,  in
differential geometry one says \textquotedblleft pull-back\textquotedblright of $p-$forms, which is 
built upon the
notion of differential of a function:

\begin{teo}
\textbf{Definitions}. Let $\phi: M\rightarrow N$ be a  map of manifolds and let
$\phi(x)=y$. Let $T_xM$ and $T_{\phi(x)}N$ be the tangent spaces at $x$ over $M$
and, respectively, at $\phi(x)$ over $N$. We define the \textbf{differential}
$\phi_* $ of $\phi$ as the linear isomorphism $\phi_* : T_xM\rightarrow
T_{\phi(x)}N$ such that for any vector $v \in T_xM$ and any scalar function $f:N\rightarrow 
\mathbb{R}$  we have $[\phi_*(v)](f) =
v(\phi\circ f)$. \textbf{A
map is
called differentiable} or smooth if its differential exists. For a smooth map  $\phi: M\rightarrow 
N$  the
\textbf{pull-back} $\phi^* : \Lambda_{\phi(x)}^1 N \rightarrow 
\Lambda_{x}^1 M$  is a linear
transformation such that for $x \in M$ we get  $\phi^*(\omega)(v) = \omega(\phi_*(v))$ for all 
vectors $v\in T_xM$ and
1-forms $\omega$.
\end{teo}

\bigskip

The pull-back observes the following properties, the first of which allows to
extend the pull-back of 1-forms to the entire space of forms:

\begin{enumerate}
\item   $\phi^*(\alpha \wedge \beta ) = \phi^*(\alpha) \wedge \phi^*(\beta )$

\item  $\phi^*(d\alpha) = d(\phi^*(\alpha) )$

\item $(\phi \circ \psi)^* = \psi^* \circ \phi^* $ 

\item $(\phi^{-1})^*=(\phi^*)^{-1}$, if $\phi^{-1}$ exists.

\end{enumerate}

\bigskip

\begin{teo}
\textbf{Lemma. } Let $\phi: M\rightarrow N$ be an orientation preserving
diffeomorphism between manifolds and let $H^*(M)$ and $H^*(N)$ the respective
cohomology rings, then $\phi$ induces a ring contravariant isomorphism

\

$F_{\phi}^*: H^*(N)\rightarrow H^*(M)$  defined by
$F_{\phi}^*([\eta])=[\phi^*(\eta)].$
\end{teo}

Proof: Let us prove that $F_{\phi}^*$ is well defined, with an inverse which is the ring 
homomorphism induced by $\phi^{-1}$ and that 

$F_{\phi}^*([\eta_1] \sqcup 
[\eta_2])=F_{\phi}^*([\eta_1]) \sqcup F_{\phi}^*([\eta_2]).$

To see that $F_{\phi}^*$ is well defined we need to prove that two
representatives of the same cohomology class in the domain are transformed into
members of the same cohomology class in the range. Two elements are in the same
cohomology class  if they differ by an exact form, i.e., if $[\eta_1]=[\eta_2]$
then $\eta_1=\eta_2+d\tau$. In this case,

$\phi^*(\eta_1)=\phi^*(\eta_2)+\phi^*(d\tau) = \phi^*(\eta_2)+d(\phi^*(\tau))$

this means that the members of a class are transformed into elements that differ by an exact form, 
i.e., they are members of the same class:

if $[\eta_1]=[\eta_2]$ then $[\phi^*(\eta_1)]=[\phi^*(\eta_2)]$ so that 
$F_{\phi}^*([\eta_1])=F_{\phi}^*([\eta_2])$ is well defined.

\bigskip

Let us verify now that if $F_{\phi}^*([\eta])=[\phi^*(\eta)]$, the inverse of $F_{\phi}^*$  is the 
ring homomorphism $G_{\phi^{-1}}^*$ induced by $\phi^{-1}$, i.e., if  
$G_{\phi^{-1}}^*([\omega])=[(\phi^{-1})^*(\omega)]$ then $G_{\phi^{-1}}^*(F_{\phi}^*([\eta])) = 
[\eta] $  and  $F_{\phi}^*(G_{\phi^{-1}}^*([\omega])) = [\omega] $. For the first part we have:

$G_{\phi^{-1}}^*(F_{\phi}^*([\eta])) = G_{\phi^{-1}}^*([\phi^*(\eta)]) = 
[(\phi^{-1})^*(\phi^*(\eta))]=
[(\phi^*)^{-1}(\phi^*(\eta))]= [I\eta]=[\eta]$. The proof of the second part is similar.

\bigskip

Let us show that $F_{\phi}^*$ is indeed a ring homomorphism:

$F_{\phi}^*([\eta_1] \sqcup [\eta_2])=F_{\phi}^*([\eta_1\wedge \eta_2]) = [\phi^*(\eta_1\wedge 
\eta_2) ]= [\phi^*(\eta_1)\wedge \phi^*(\eta_2) ]$

$=[\phi^*(\eta_1)]\sqcup[ \phi^*(\eta_2) ]$
$=F_{\phi}^*([\eta_1]) \sqcup F_{\phi}^*([\eta_2])$.

In conclusion, $F_{\phi}^*$ is a ring homomorphism:

$F_{\phi}^*([\eta_1] \sqcup [\eta_2])=F_{\phi}^*([\eta_1]) \sqcup F_{\phi}^*([\eta_2])$.

\bigskip

We also say that $F_{\phi}^*$ is orientation preserving, in the sense that if we
choose a basis $E$ of $H^{2k}(N)$ then $F_{\phi}^*(E)$ will be a basis of
$H^{2k}(M)$ that observes the same orientation.

\begin{teo}
 \textbf{Convention. } The ring isomorphism $F_{\phi}^*$ is denoted as $\phi^*$.
\end{teo}

\begin{teo}
\textbf{Theorem. } The  quadratic signature  form $Q$ is invariant under diffeomorphisms.
\end{teo}

Proof.  Let $\phi: M\rightarrow N$ be a diffeomorphism between manifolds, and let $F_{\phi}^*$ be 
its induced isomorphism between $H^*(M)$ and $H^*(N)$. The theorem of change of variables reads:

$\int_{N=\phi(M)} \eta = \int_M \phi^*(\eta)$.

Let us apply this theorem to the signature quadratic form:

$\int_{N=\phi(M)} [\eta_1]\sqcup [\eta_2] = \int_{N=\phi(M)} \eta_1 \wedge \eta_2 =\int_M  \phi^* 
(\eta_1 \wedge \eta_2 )= \int_M  \phi^*\eta_1 \wedge \phi^*\eta_2 = $

$ \int_M  [\phi^*\eta_1 \wedge \phi^*\eta_2]) = \int_M  [\phi^*\eta_1])\sqcup ([\phi^*\eta_2]) =$

$\int_M F_{\phi}^*( [\eta_1])\sqcup F_{\phi}^*([\eta_2])$.

Thus, we have proved that $\int_{N=\phi(M)} [\eta_1]\sqcup [\eta_2] = \int_M F_{\phi}^*( 
[\eta_1])\sqcup F_{\phi}^*([\eta_2])$ and this means that two diffeomorphic manifolds have the same 
quadratic signature form: 

\

$ Q([ \eta_1],  [\eta_2]) = \int_{N=\phi(M)} [\eta_1]\sqcup [\eta_2] = \int_M F_{\phi}^*( 
[\eta_1])\sqcup F_{\phi}^*([\eta_2]) = Q([\phi^* \eta_1],  [\phi^* \eta_2])$.

\bigskip

\begin{teo}
\textbf{Corollary. } Diffeomorphic oriented manifolds  have  the same signature.
\end{teo}

The following illustration explains the whole idea. Let us suppose, 
as in an example above,  that the signature quadratic form in $N$ with respect to basis $\{ 
d\theta, d\phi \}$ has matrix

$$Q =
\bordermatrix{~ &d\theta  & d \phi  \cr
                  d\theta & 0 & c \cr
                  d \phi  & -c & 0\cr}$$
                  
and that we have a diffeomorphism  $\phi: M\rightarrow N$   that 
preserves orientation. Now, the preimage of any basis $E$ of $H^{2k}(N)$ through $\phi^*$ is also a 
basis, otherwise 
$\phi$ would not be a diffeomorphism. Therefore, the matrix of the signature form in $M$ with 
respect to basis        $\{\phi^* (d\theta), \phi^*(d\phi) \}$ has matrix

$$Q =
\bordermatrix{~ &\phi^* (d\theta)  & \phi^*(d\phi)   \cr
                  \phi^* (d\theta) & 0 & c \cr
                  \phi^*(d\phi)   & -c & 0\cr}$$           

We see that the entries of the matrix are conserved. In particular,  signs are  conserved  as a 
result of the orientation preserving property of  $\phi$. Changes affect
only the labels of the matrix and therefore the structure of eigenvalues and eigenvectors is 
untouched as it is also its signature.

\subsection{Invariance of $\sigma(M)$ under homotopy equivalence}

Conceptually, smooth manifolds are inseparable from smooth functions: What
does happen with the signature   if the manifold is smoothly deformed?

\bigskip

\begin{teo}
\textbf{Definition. } Two functions with the same domain and codomain, $h_0,h_1:
Z\rightarrow W$, are \textbf{homotopically  equivalent}, $h_0\sim h_1$, if there
exists a smooth map
$F: Z\times [0,1] \rightarrow W $ such that  $F|_{Z\times \{0\}} =h_0$
and $ F|_{Z\times
\{1\}} =h_1$.
\end{teo}

 Remarks:  Homotopic equivalence is  a topological
concept so, the closed interval $[0,1]$ is endowed with the relative  or subspace 
topology inherited from   $(-\epsilon, 1+ \epsilon)$.  Intuitively , two functions are
homotopically equivalent if there
exists a continuous deformation of $h_0(Z)$ into $h_1(Z)$ or just imagine
yourself taking a bar of plasticine and
deforming it continuously from an initial state into a final one. Now,
we changed continuity for smoothness because we deal with smooth manifolds and
so we demand from F
to be also smooth. 
The functions shall not be onto and the images of the two functions could be
disjoint. 

\bigskip

\begin{teo}
\textbf{Definition. } Let $M$ and $N$ be two oriented smooth compact manifolds. 
We say that $M$ and $N$ are \textbf{(strongly) homotopically equivalent} if
there exist two orientation preserving smooth maps

$$g:N\rightarrow M$$

$$f:M\rightarrow N$$

such that $f\circ g \sim id_N$ and $g\circ f \sim id_M$.
\end{teo}

To understand the meaning of this concept let us imagine that $M$ and $N$
represent plasticine figures.   If $M$ can deform itself
smoothly into its image $g(f(M))$ and if $N$ also can do the equivalent in its
side, then we say that the two manifolds are homotopically equivalent.

\begin{teo}
\textbf{Question. }   Let us consider the following two manifolds, the first a
projective space, the second a torus: $ \mathbb{CP}^3 $  and $S^2\times S^4$.
These two manifolds seem very similar
according to certain descriptors. Let us see. 
\end{teo}

To begin with, both have  real dimension 6. Moreover, they both have the same
associated polynomial:

From our construction of  $ \mathbb{CP}^n$, we have the following decomposition:

\

 $ \mathbb{CP}^3  =  e^0 \cup e^2  \cup e^{4} \cup e^6 \leftrightarrow p(t)  = 1
+ t^2 + t^4 + t^6$.

\

On the other hand, the decomposition of the super-torus $S^2 \times S^4$ can be
calculated thanks to Poincaré polynomials:

\

$ S^2 \times S^4   \leftrightarrow 
(1+t^2)(1+t^4) = 1 + t^2 + t^4 + t^6   \leftrightarrow  e^0 \cup e^2  \cup e^{4} \cup e^6$.

\

We see that these two spaces have the same Poincaré polynomial, so they share  
the
same cohomology.     Does this means that they are homotopy  equivalent? To
answer this question,  we will show that homotopy
deformations  generate isomorphic cohomology rings.
As a consequence,   the signature of a manifold is conserved under homotopy
equivalence. In conclusion, our two manifolds will be not  homotopy equivalent
if they have non isomorphic cohomology rings. 

\

To begin with, let us prove that the antitransport of closed forms through two
homotopy equivalent functions belong in the same cohomology class, i.e.,   that
they differ by an exact form.

\begin{teo}
\textbf{Lemma. } Let $f,g: M\rightarrow N$ be  smooth maps  that are homotopic to each other. If 
$\omega\in \Omega^K(N)$ is a closed form, the difference of the  pull-back images is exact:

$$f^*\omega -g^*\omega = d\psi$$

where $\psi \in \Omega^{k-1}(M)$ and  $f^*$ and $ g^*$ are the pull-backs of
$f$ and $g$ respectively.
\end{teo}

Proof. Since $f\sim g$, there exists a smooth map
$F: M\times [0,1] \rightarrow N $ such that

$F|_{M\times \{0\}} =f$  and$ F|_{M\times \{1\}} =g$, i.e., $F(x,0)=f(x)$
and $F(x,1)=g(x)$ for $x\in M$.

\bigskip

Now we will be involved in a game that uses $F$ and the fundamental theorem of calculus: to prove 
that $f^*\omega -g^*\omega = d\psi$ for some $\psi$, we will prove a rigorous version of the 
following idea: $f^*\omega -g^*\omega = d\int F^*\omega$. In what follows in this lemma, the 
integral will be replaced by the operator $P$.

\bigskip

Let us consider a $k-$form $\eta \in \Omega^k(M\times [0,1])$. $\eta$ takes the form

$\eta = a_{i_1...i_k}(x,t)dx^{i_1}\wedge...dx^{i_k} + b_{j_1...j_{k-1}}(x,t)dt\wedge 
dx^{j_1}\wedge...\wedge dx^{j_{k-1}}$

where $x\in M, t\in [0,1]$. The second term is of degree $k$ but it shall include $dt$ so it has 
only $k-1$ degrees of freedom to choose its components.

\bigskip

Define a map $P:\Omega^k(M\times I)\rightarrow \Omega^{k-1}(M)$ by

 $P ({\eta }) = (\int_0^1  b_{j_1...j_{k-1}}(x,s)ds) dx^{j_1}\wedge...\wedge dx^{j_{k-1}}$

Let $f_t$ be a map $f_t:M\rightarrow M\times I$ such that $f_t(p)=(p,t)$

We have

$f_t^*\eta = a_{i_1...i_k}(x,t)dx^{i_1}\wedge...dx^{i_k} \in \Omega^k(M)$

since $f_t^*(dt\wedge dx^{j_1}\wedge...\wedge dx^{j_{k-1}})=0$.

Let us prove now that

$d(P ({\eta})) + P(d\eta) = f_1^*(\eta)-f_0^*\eta$.

Indeed, if we calculate each term in the lhs we get

$d(P({\eta}))= d(\int_0^1  b_{j_1...j_{k-1}}(x,s)ds) dx^{j_1}\wedge...\wedge dx^{j_{k-1}}$

$=(\int_0^1  (\partial b_{j_1...j_{k-1}}(x,s)/\partial x^{j_{k}})ds) dx^{j_{k}}\wedge 
dx^{j_1}\wedge...\wedge dx^{j_{k-1}}$

\bigskip

On the other hand

$P(d\eta) = P[d(a_{i_1...i_k}(x,t)dx^{i_1}\wedge...dx^{i_k} + b_{j_1...j_{k-1}}(x,t)dt\wedge 
dx^{j_1}\wedge...\wedge dx^{j_{k-1}}]$

$=P[(\partial a_{i_1...i_k}(x,t)/\partial x^{i_{k+1}})dx^{i_{k+1}}\wedge dx^{i_1}\wedge...dx^{i_k}$

$+ (\partial a_{i_1...i_k}(x,t)/\partial t)dt\wedge dx^{i_1}\wedge...dx^{i_k} ]$

$+(\partial b_{j_1...j_{k-1}}(x,t)/\partial x^{j_{k}})dx^{j_k}\wedge dt\wedge 
dx^{j_1}\wedge...\wedge dx^{j_{k-1}})$

$=(\int_0^1 (\partial a_{i_1...i_k}(x,s)/\partial s)ds)dx^{i_1}\wedge...dx^{i_k}$

$  -(\int_0^1 (\partial b_{j_1...j_{k-1}}(x,s)/\partial x^{j_{k}})ds) dx^{j_k}\wedge 
dx^{j_1}\wedge...\wedge dx^{j_{k-1}}$

\bigskip

Summing up these two terms we get

$d(P({\eta})) + P(d\eta) =  (\int_0^1 (\partial a_{i_1...i_k}(x,s)/\partial 
s)ds)dx^{i_1}\wedge...dx^{i_k}$.

Applying the fundamental theorem of calculus for a continuous function
we have

$d(P({\eta})) + P(d\eta) =[a_{i_1...i_k}(x,1)-a_{i_1...i_k}(x,0)]dx^{i_1}\wedge...dx^{i_k}$

$d(P({\eta})) + P(d\eta) = f_1^*(\eta)-f_0^*\eta$

Let us apply this identity to the pull-back of a closed form  $\omega \in \Omega^k(N)  $ and

$\eta = F^*\omega \in \Omega^k(M\times[0,1])$:

so $d(P({\eta})) + P(d\eta)$ becomes

$d(P({F^*\omega})) + P(dF^*\omega) = f_1^*(F^*\omega)-f_0^* (F^*\omega) $

recalling that $(FG)^* = G^*F^*$ we get

$d(P({F^*\omega})) + P(dF^*\omega) = (Ff_1)^*\omega-(Ff_0)^*\omega $

but $f_t(x) = (x,0)$ and moreover $f\sim g$ so that $F(f_0(x))=F(x,0)=f(x)$ and 
$F(f_1(x))=F(x,1)=g(x)$

Hence

$d(P({F^*\omega})) + P(dF^*\omega) = f^*\omega -g^*\omega$.

Recalling now that $\omega$ was chosen to  be closed, $d\omega = 0$, we get 
$F^*d\omega =0$. We can rewrite this as $d(F^*\omega) =0$ because the pull-back
and the exterior derivative commute.  Integrating with $P$ between 0 and 1 we
have  $P(dF^*\omega) = 0$. Replacing

$d(P({F^*\omega})) + P(dF^*\omega) = d(P({F^*\omega})) = f^*\omega -g^*\omega$

\bigskip

Explicitly $ f^*\omega -g^*\omega = d\psi $ where $ \psi = P({F^*\omega})$ proving that the 
difference of the pull-backs of closed forms throughout homotopy equivalent functions is exact.

\begin{teo}
\textbf{Corollary. } Let $f,g: M\rightarrow N$ be maps which are homotopic to
each other. Then, the pull-back maps $f^*, g^*:H^*(N) \rightarrow H^*(M)$
defined on the de Rham cohomology rings are identical, i.e., for a closed form $\omega\in
\Omega^K(N)$, we have  $[f^*\omega] = [g^*\omega ] $.
 
\end{teo}

In fact,

$  [f^*\omega] -[g^*\omega ] = [f^*\omega -g^*\omega ]=[d\psi]=0 $

because the zero class of $H^*(M)$ is conformed by    the forms that are exact.

\begin{teo}
\textbf{Theorem. } If $M$ and $N$ are homotopy equivalent,   oriented, even
dimensional, compact manifolds, then $\sigma(M) = \sigma(N)$.
\end{teo}

Proof. Since $M$ and $N$ are homotopy equivalent manifolds, there exist two orientation preserving 
smooth maps

$f:M\rightarrow N$

$g:N\rightarrow M$

such that $f\circ g \sim id_N$ and $g\circ f \sim id_M$. We can now apply the
result of the previous corollary: when two functions are homotopic to each
other, their pullbacks are identical. So, on one hand we have:

 $(f\circ g)^* =  g^*\circ f^* = (id_N)^* = id_{H^*(N)}$ 
 
 and on the other

 $(g\circ f)^* = f^*\circ g^*  = (id_M)^* = id_{H^*(M)}$.

 \
 
 These system of equalities says us that as $f^*$ as $g^*$ are isomorphisms or
that $M$ and $N$ are diffeomorphic. Since they are oriented, their ring
structure is also isomorphic and henceforth they have the same signature.

\

Next definitions and ensuing comment teach us how to simplify manifolds to
its most fundamental  cores without distorting their differential structure.
Say,  a torus $S^1  \times S^1$ can be considered as the simplification of all
those surfaces into which it can be smoothly deformed. Compare with Zeeman
(\cite{Zeeman66}, 1966).

\begin{teo}
\textbf{Definition. } Let   $R$ be a, not empty, topological subspace of  $M$. If
there exists a continuous map $f:M\rightarrow R$ such that $f|_{R} = id_R$, $R$
is called a \textbf{retract} of $M$ and  $f$ a \textbf{retraction}.
\end{teo}

A retraction is our formalization of a curvilinear projection. 

\begin{teo}
\textbf{Definition. } Let   $R$ be a, not empty, topological subspace of  $M$ and
$f$ a retraction of $M$ over $R$. $R$ is said to be a \textbf{deformation
retract} if
$id_M$ and $f$ are homotopically equivalent and $R$ is point by point invariant
in the deformation.
\end{teo}

We shall highlight the fact that our retractions
eliminate homotopic redundances but holes are not eliminated. So, a circle
cannot be joined to its north pole by a deformation retraction.

\bigskip

\textbf{Disclaimer}. Two manifolds could be inequivalent and yet they can have
the same cohomology vector spaces and henceforth the same signature. Example:
take as
$M= S^2\times S^4$ and as $N = \mathbb{CP}^3$. To prove that they are
homotopiclly inequivalent, we will show that they do not  expand isomorphic
cohomology
rings.   In fact,  they are dissimilar:  $ \mathbb{CP}^3$ is recursive up to 3
complex
dimensions, i.e., 6 real ones,   while $ S^2 \times S^4 $ is recursive only up
to 4
real dimensions. Thus products vanish above 4 terms in the last case with the
exception of
that product that corresponds to the volume form, while we can rise up to 6 in
the
former one.

\begin{teo}
\textbf{Summary}. The maximal subspaces of forms over which $Q$ is
definite-positive have the same dimension for oriented, boundariless, compact,  
 finite dimensional = 4k, homotopically equivalent manifolds.
\end{teo}

\subsection{The Hodge star operator}

In the sequel we will present an elaboration of a nice observation regarding
our matrices $Q$: the only products that matter are those whose output
complete the volume form. Thus, a question arises: can we   define an operator
that associates to any given form what it lacks to complete the volume form?
The solution to this problem has shown to be very rich 
if we consider a closed riemannian manifold with a metric $(M,g)$, where we can define the Hodge-* 
operator as above. We need some few properties of this operator.

\begin{teo}
\textbf{Properties of *}:
\end{teo}

We  need the following  fundamental  properties of the Hodge-* operator,
where $vol$ is the volume form
(Dray, \cite{Dray99} 1999;  Ivancevic et al, \cite{Ivancevic11} 2011):

\begin{enumerate}
\item $\alpha \wedge *\beta = (\alpha, \beta) vol$.

 \item The star operator provides what  a form lacks to be the volume form:
 
 $\alpha \wedge *\alpha = ||\alpha  ||^2 vol$
 
 where   $ ||\alpha  ||^2 = (\alpha, \alpha)$
 
 \item $** =  (-1)^{p(n-p)}$

\

and hence

\

$**(-1)^{p(n-p)} = 1$

\

This implies that

\

 $$ * ^{-1} = (-1)^{p(n-p)} *$$ 

\

\item We also need $*^*$, the adjoint of * for the scalar product between
p-forms. That product  reads:

\

 $(\alpha , \beta ) = \int_\Omega \alpha \wedge * \beta$.

\

 From the property

\

$(*\alpha, *\beta ) = (\alpha, \beta) $

\

we get

\

$(*\alpha, *\beta ) = (\alpha,*^* * \beta)  = (\alpha, \beta) $

\

Hence $*^* * = I$. This implies that

\

 $$ * ^{-1} = (-1)^{p(n-p)} * = *^*$$ 

\

\begin{teo}
\textbf{Example. } Let us calculate $*$ over some forms of $\mathbb{R}^{4}$.
\end{teo}

To calculate $*$, we use the next trick: $(*\alpha)$ must complement $\alpha$ to
fill in  the volume form $dxdydzdv$ and the sign must be adjusted accordingly.
Examples:

\

$*dxdy = dzdv$ because $(dxdy)(dzdv) =  dxdydzdv$.

$*dxdzdv = dy$ because $ dxdzdvdy = -dxdzdydv = dxdydzdv$.

$*dydzdv = -dx$  because

$dydzdv(-dx) = -dydzdvdx = dydzdxdv = -dydxdzdv= dxdydzdv$.

\

 \begin{teo}
\textbf{Example. } Let us inquire in $\mathbb{R}^{4}$ over the eigenvectors of
$*$.
\end{teo}

Since $*$ completes forms to fill in the volume form, possible eigenvectors of
$*$ must be a linear combination of terms that complete one another. So, let us
prove that $\omega  = dxdy +  dzdv$ is an eigenvector of $*$.
In fact:

\

$*\omega = *(dxdy +  dzdv) = dzdv + dxdy = \omega$

\

Thus, $\omega$ is an eigenvector with eigenvalue 1. To fabricate an eigenvector
with eigenvalue -1, let us try $\eta = dxdz +  dydv$:

\

$*\eta = *( dxdz +  dydv) = - dydv-dxdz = -\eta.$

\

\subsection{$\delta $: the adjoint of the derivative $d$}

We define $\delta = d^* $, as the adjoint of $d$ which is defined by
the equation

\

$(d\alpha, \beta) = (\alpha, \delta \beta)$.

\

It is found that $\delta = (-1)^{n(p+1)+1} *d*$. When $n$ is even, $\delta = -
*d*$

\begin{teo}
\textbf{Theorem. } $\delta^* = d$ (for $n$ even).
\end{teo}

To prove this theorem, we begin with

\

$\delta =  - *d*$

\

to get

\

$\delta^*  =  -  (*d*)^* =  - *^* \delta *^* = - (-1)^{p(n-p)}
(-1)^{p(n-p)}*\delta *$
$=   **d** =  d $.

\end{enumerate}

\

\

\subsection{The Euler characteristic}

Because $** = (-1)^{k(n-k)}$, when operating over $k-$forms, there is a natural  bijection between 
$H^k$ and $ H^{n-k}$ known as Poincar\'e' duality. Example $dx $ and $dydz$ single out one to 
another in $\mathbb{R}^{3}$. Let the \textbf{k-Betti number} $b_k$ be defined by 
 $b_k = dim H^k$. Define the \textbf{Euler characteristic} of $M$ as $\chi(M) = \sum (-1)^kb_k$.

\begin{teo}
\textbf{Lemma. } Let $dim M= 4k$ then the Euler characteristic of $M$ and the dimension of $H^{2k}$ 
have the same parity, i.e., $\chi(M)=b_{2k}$ mod 2.
\end{teo}
Proof. $\chi(M) = \sum_0^{4k} (-1)^kb_k =\sum_0^{2k-1} (-1)^kb_k+ (-1)^{2k}b_{2k}+ \sum_{2k+1}^{4k} 
(-1)^kb_k$.

Using the Poincar\'e duality, this can be rewritten as

 $\chi(M) = 2\sum_0^{2k-1} (-1)^kb_k + (-1)^{2k}b_{2k}$

$\chi(M)-b_{2k} = 2\sum_0^{2k-1} (-1)^kb_k  $

that shows that at both sides of this equation we are dealing with even numbers. Or, equivalently, 
$\chi(M)=b_{2k}$ mod 2.

\section{Harmonic forms}    

It is shown here that there is a suitable basis for  the calculation 
of the signature of a manifold.

\begin{teo}
\textbf{Definitions}. The \textbf{Dirac operator} $D$ is $D= d+\delta$ and the
\textbf{Laplacian} $\triangle = D^2 = (d+\delta)^2 = (d+\delta)(d+\delta) = d^2
+ d\delta + \delta d +\delta^2 = d\delta + \delta d$. A form that satisfies the
\textbf{Laplace equation} $\triangle \omega =0$ is called \textbf{harmonic}. The
space of harmonic forms of degree $k$  is denoted as $Harm^k(M)$.
\end{teo}

\begin{teo}
\textbf{Example. } Let us exhibit the harmonic representatives of $H^k$ over
$S¹$.
\end{teo}

For $S^1$, $b_o = 1$, since it has just one connected component, and  $b_1 = 1$
since it has just one 1-dimensional hole. We will look for 
solutions to the equation $\nabla \alpha = (d\delta + \delta d )(\alpha)$.
Because $S^1$ has dimension 1, we need the basic definition of $\delta =
(-1)^{n(p+1)+1} *d*$ that for $n= 1$ becomes $\delta = (-1)^{p+2} *d* = 
(-1)^{p}*d*$. So, we are looking for solutions of 

\

$\nabla \alpha = (-1)^{p}(d\delta + \delta d )(\alpha) =  (-1)^{p}(d*d*
+*d*d)(\alpha) = 0$.

\

or of

\

$(d*d* +*d*d)(\alpha) = 0$.

\

The obvious candidates for harmonic forms are: the constant function 1, a
0-form,  and $d\theta$, a 1-form. Let us inquire whether or not  they are
harmonic.
We shall use the fact that $*0 = 0$ because $*0 = 0d\theta = 0$.

\

Let us test the 0-form $1$:

\

$(d*d* +*d*d) 1 = d*d*1  + d*d1 =  d*dd \theta +*d*0 = d*0 +0 = 0$

\

So, $1$ is harmonic. What happens with the 1-form $1d\theta$?

\

$(d*d* +*d*d) (1d\theta) = d*d*(1d\theta) +*d*d(1d\theta) =  d*d1 +0 = 0 $

\

So, $d\theta$ is harmonic. 

\begin{teo}
\textbf{Example. } Let us exhibit the harmonic representatives of $H^k$ over
$S^2$. We consider the usual coordinates $\theta$ and $\phi$ in that order.
\end{teo}

$S^2$ has one connected component, so $b_o=1$. It has one 2-hole, so $b_2 = 1$
and every closed curve over it is contractile to a point so, $b_1 = 0$. In
consequence, let us show that the constant function $1$ and the volume form
$d\theta d\phi$ are harmonic. With $1$ we have:

\

$(-d*d* -*d*d) 1 = -d*d*1 -*d*d1 =  -d*dd \theta d\phi -*d*0 = -d*0 -0 = 0$

\

So, $1$ is harmonic.  What happens with $1d\theta d\phi$?

\

$(-d*d* -*d*d) (1d\theta d\phi) = -d*d*(1d\theta d\phi)  -*d*d( 1d\theta d\phi )
=  -d*d1 - 0 = 0 $

\

So, $d\theta d\phi$ is harmonic.

\

\begin{teo}
\textbf{Theorem. } $D^* = D$.
\end{teo}

We have:

\

$D = d +\delta $

\

so

\

$D^*  = d^* + \delta^*  = \delta + d = D$.

\begin{teo}
\textbf{Theorem. } $\triangle$ is  selfadjoint.
\end{teo}

Proof: $(\triangle \omega, \psi) = ( D^2\omega, \psi)=$
$( D\omega, D^*\psi) = ( D\omega, D\psi) = ( \omega, D^* D\psi) = ( \omega, D^2\psi) =   ( \omega, 
\triangle\psi)$.

\begin{teo}
\textbf{Remark. } We have proved that $(\triangle \omega, \psi) = ( \omega,
\triangle\psi)$ but only if that  makes sense.   So, we  need   to give  a
glance at the domain of definition of involved operators. Basically, we have two
items: integration and derivation  of differential forms. Since at last these are
reduced to integration of functions involving       partial derivatives of real
valued functions that are defined  over open sets  of  $\mathbb{R}^n$,  
domains must be referred to them. Moreover, our arguments strongly rely on duality so, we must 
resort to Sobolev spaces.  In       the       case of the Dirac operator we must think of
Sobolev space  of order 1, and  of 
Sobolev space  of order 2 for the   Laplacian. Sobolev spaces are complete so, the usual inner 
product makes then  into 
Hilbert Spaces (Paycha, \cite{Paycha97} 1997).  Now, to define Sobolev spaces intrinsically, we 
must consider patching local definitions thanks to partitions of unity.   Actually, this trickery
is already present in the integration machinery.                          
                                                                               
\end{teo}

\begin{teo}
\textbf{Theorem. } $\triangle$ is positive. In other words: $(\triangle \omega, \omega) \geq 0 $ 
where $\omega\in \Omega(M)$.
\end{teo}

Proof. $(\triangle \omega, \omega) = ((d\delta+ \delta d) \omega, \omega)$
=$(d\delta \omega, \omega)+ (\delta d \omega, \omega) = (\delta \omega, \delta \omega)+ (d \omega, 
d\omega) = \|\delta \omega \|^2 + \|d \omega \|^2  \geq 0$

\begin{teo}
\textbf{Corollary. } $(\triangle \omega, \omega) = 0$ iff $\delta \omega= d \omega=0$, i.e., 
harmonic forms are both closed and coclosed.
\end{teo}

\begin{teo}
\textbf{Theorem.} Poisson's equation $\triangle \psi = \omega $ has a solution iff $\omega $ is 
orthogonal to $Harm(M)$. The solution $\psi$ is noted as $\psi = \triangle^{-1}\omega$.
\end{teo}

Proof.  Let us assume that there exists $\psi $ such that  $\triangle \psi = 
\omega $. Then, 
by taking $\gamma$ from $Harm(M)$ we get:

$(\omega, \gamma) = (\triangle \psi, \gamma) = ( \psi, \triangle\gamma)=(\psi, 0) =0$

which shows that the orthogonality of $\omega$ to $Harm(M)$ is necessary for the existence of a 
solution of the Poisson's equation . On the other hand, if $\omega  $ is 
orthogonal to $Harm(M)$, its preimage cannot be in $Harm(M)$ because it is the kernel of 
$\triangle$. Hence,  every preimage must have a component that belongs in $Harm(M) ^ \perp $, the 
orthogonal complement of $Harm(M)$. Let us consider the restriction of $\triangle $ to $Harm(M) 
^ \perp $, i.e., tolerating an abuse of notation, we consider $\triangle: Harm(M) ^ \perp  
\rightarrow Harm(M) ^ \perp $: let us prove that this operator is one to one and onto. It is in this 
sense that $\triangle$ is 
a bijection with an inverse that can be denoted as $\triangle^{-1}$.

To see why $\triangle: Harm(M) ^ \perp  \rightarrow Harm(M) ^ \perp $ is one to one, let us take  
$\phi$ and $  \psi$  in $Harm(M) ^ \perp $. In 
this case, $ \phi -  \psi$ is also in  $Harm(M) ^ \perp $. Suppose now that $\triangle \phi = 
\triangle \psi$. Then, $\triangle (\phi -  \psi) = 0 $. So, $\phi -  \psi$ is in $Harm(M)$. Now, the 
only element that is in both $Harm(M)$ and $Harm(M) ^ \perp $ is the zero element. Henceforth, $\phi 
=  \psi$.

To prove that $\triangle: Harm(M) ^ \perp  \rightarrow Harm(M) ^ \perp $ is onto and that 
therefore $\triangle \psi = \omega $ has a solution, we 
apply the Riesz' 
Representation Theorem for continuous linear functionals on Hilbert spaces. In fact, the equation 
$\triangle \psi = \omega $  implies

$(\triangle \psi, \phi) = (\psi, \triangle \phi) =  (\omega, \phi) $  for $\phi $.

Now, because $\omega $ is a   fixed element, the expression $(\omega, \phi)$ defines a linear 
functional that can be denoted as $l_{\omega} (\phi) = (\omega, \phi)$. So, we can define over 
$Harm(M) ^ \perp $ the inner product $[[, ]]$  given by:

$[[\psi, \phi]] = ( \psi, \triangle \phi)$

The Riesz' 
Representation Theorem allows us to represent the functional $l_{\omega} (\phi) = (\omega, \phi)$ by 
a single element $\Omega$ that 
operates through   the inner product $[[,]]$: 

$(\omega, \phi) = l_{\omega}(\phi) =  [[\Omega, \phi]] = (\Omega, \triangle\phi)  = 
(\triangle\Omega, 
\phi)  $ 

Since this is certain for every $\phi$, we conclude that  $ \triangle\Omega = \omega$ and that 
therefore the Poisson's equation has a solution when $\omega$ is orthogonal to $Harm(M)$. Now, a 
lot of analysis over Sobolev spaces is necessary to fill in the details in regard with continuity, 
see  (Min Ru, \cite{MinRu00}, 2000).

\begin{teo}
\textbf{Theorem. } Harmonic, Exact and co-exact forms are mutually orthogonal spaces.
\end{teo}

Proof. Since $d^2=0$ then $(d^2 \alpha_{k-1}, \beta_{k+1}) =(d \alpha_{k-1}, \delta \beta_{k+1}) 
=0$, showing that exact and coexact forms are orthogonal.

Since harmonic forms are closed then if $\gamma_k $ is harmonic then

$d\gamma_k=0$ and $(\beta_{k+1},d\gamma_k) =(\delta \beta_{k+1},\gamma_k)= 0$ showing the mutual 
orthogonality between coexact and harmonic forms.

Likewise, harmonic forms are closed, i.e.,   $\delta\gamma_k=0$. Then $(\alpha_{k-1},\delta 
\gamma_k) =(d \alpha_{k-1},\gamma_k)= 0$ showing the mutual orthogonality between exact and 
harmonic forms.

\begin{teo}
\textbf{Hodge Decomposition Theorem}. Harmonic, Exact and co-exact forms generate $\Omega^k(M)$: 
$\Omega(M)=Ker\triangle \oplus Ker d \oplus Kerd^* $.
\end{teo}

Proof. Let $P:\Omega^k(M)\rightarrow Harm^k(M)$ be the orthogonal projection operator generated by 
the scalar product of forms, then for any $\omega \in \Omega^k(M)$ we have that $\omega -P\omega$ 
is orthogonal to $Harm^k(M)$. Hence,  Poisson's equation

$\triangle \psi =  \omega -P\omega$ has a solution that can be written as $ \psi =  
\triangle^{-1}(\omega -P\omega)$

In conclusion, it makes sense to write

$\omega = \triangle \psi  + P\omega =(d\delta + \delta d)\psi  + P\omega  = d(\delta \psi ) + 
\delta (d\psi)  + P\omega$

that reads:  any form can be orthogonally decomposed as a sum of an exact form plus a coexact form 
plus a harmonic form.

\begin{teo}
\textbf{Hodge  Theorem}. $H^k(M) = Harm^k(M)$, in other words, any cohomology class is the class of 
a harmonic form and every harmonic form is a nontrivial member of $H^k(M)$.
\end{teo}
Proof.
Our universe is the space of closed forms, since $H^k(M)$ is the space of equivalence classes of 
closed forms defined by the relation $\omega \sim \omega'$ if there exists $d\psi $ such that 
$\omega = \omega' + d\psi$.

Let us show that there is a isomorphism between $H^k(M)$ and $ Harm^k(M)$ induced by the projection 
operator over harmonic forms $P: \Omega^k(M)\rightarrow Harm^k(M)$. The isomorphism is

$P([\omega]) = [P\omega]$.

Decomposing $\omega $ in its components:

$\omega= d(\delta \psi ) + \delta (d\psi)  + P\omega$

but $d\omega = 0$ hence applying $d$ to both sides of this equation

$d\omega= d^2(\delta \psi ) + d\delta (d\psi)  + dP\omega=d\delta (d\psi)  + dP\omega $

but $dP\omega=0 $ because $P\omega $ is harmonic and harmonic forms are closed, so

$d\omega= d\delta (d\psi) =0$

hence $(d\delta (d\psi), d\psi)  = (\delta (d\psi), \delta (d\psi))=0$ hence $\delta d\psi=0$.

Therefore

$\omega= d(\delta \psi ) + P\omega$.

Thus,  any closed form can be decomposed into an exact and a harmonic form. Moreover, we can 
clearly specify $\psi$: since $\omega- P\omega= d(\delta \psi ) $ is exact, it is orthogonal to 
harmonic forms, and so the corresponding Poisson's equation has a solution. This decomposition can 
be translated into the language of cohomology classes:

The closed form $\omega$ is in its class $[\omega]$ and another representative reads $\omega'= 
\omega + d\phi$. Proceeding as before, we find

$\omega + d\phi= d(\delta \psi' ) + P(\omega+d\phi)$

But $d\phi$ is exact, so it has no harmonic component:

$P(\omega+d\phi)=P\omega+Pd\phi= P(\omega)$

which means that if two closed forms belong in the same cohomology class, they have the same 
harmonic component. Therefore, the projector is well defined over cohomology classes.

\bigskip

Let us now  prove that different cohomology classes correspond with  different harmonic forms.

Let us suppose that $P([\omega]) = P([\eta])$. Let us decompose both forms:

$\omega= d(\delta \psi) + P\omega$

$\eta= d(\delta \rho) + P\eta$

then

$\omega - \eta = d(\delta \psi) + P\omega - d(\delta \rho) - P\eta = d(\delta \psi)  - d(\delta 
\rho)  = d(\delta \psi  - \delta \rho) $

and we have discovered that if two forms have the same harmonic component, then they are related by 
an exact form, and so they belong in  the same cohomology class. Nice.

Let us show now that any harmonic form is a nontrivial member of $H^k(M)$ that is, for $\gamma \in 
Harm^k(M)$ we have $d\gamma=0$ and $\not \exists \psi$ such that $\gamma =d\psi$. The first 
requirement is automatically fulfilled because any harmonic form is closed and the second relies on 
the fact that harmonic forms and exact forms are orthogonal one to another.

\begin{teo}
\textbf{Theorem. } $Dim Harm^k(M) = b_k < \infty$.
\end{teo}

The space of cohomology classes  $H^{k}=Ker(d_{k+1})/Im(d_k)$ is a vector space that is equal to 
 $ Harm^k(M)$, the space of harmonic forms. So, both have the same dimension.  This 
implies that if one of them has finite dimension then the other also. That of $H^k(M)$ 
is $b_k$, which  counts holes and components. Now,  $b_k$ is finite  
because   a compact manifold cannot have an infinite number of holes  or an 
infinite number of components: otherwise it would be possible to construct in both cases an open 
covering that has no finite subcovering.  Another proof, based on analysis,  that  $H^k(M)$ 
has finite dimension can be 
found in Michor(2008, pag 153, \cite{Michor08}).

\begin{teo}
\textbf{Commutation properties for $*,d,\delta$}
\end{teo}

We can decompose $\delta$  as  $\delta = (-1)^{n(p+1)+1} *d*$. In our case, $n$ is even, so   
$\delta = -*d*$ . Observe that $\delta$ takes $k-$forms and produces $k-1-$ forms. In this respect, 
it is similar to integration. Indeed, over $p-$forms $*$ produces an $(n-p)$-form over which $d$ 
produces a $(n-p+1)$-form and hence $\delta = -*d*$ produces a $n-(n-p+1)$-form or a $p-1$ form. 
Henceforth, recalling that over $p$-forms $** =(-1)^{(n-p)p}$   we have that over $(p-1)$-forms it 
reads $** =(-1)^{(n-p+1)(p-1)}$. Hence

\

$*\delta = -**d* = (-1)^{(n-p+1)(p-1)+1}d*$

\

 On the other hand, $\delta* = -*d** = -*d(-1)^{p(n-p)}=(-1)^{p(n-p)+1}*d$. Or, $\delta* = 
(-1)^{p(n-p)+1}*d$. Therefore, multiplying by the sign in the right side, we get

\

$*d= (-1)^{p(n-p)+1}\delta *$.

\

With the help of these identities we can prove the expected result:

\begin{teo}
\textbf{Theorem. } The wedge product of harmonic forms is harmonic and the Hodge star of a harmonic 
form is also harmonic.
\end{teo}

Proof: Let us take two harmonic forms, $\gamma_1, \gamma_2$. Both are closed and coclosed, i.e.,
$d\gamma_1 = d \gamma_2 = \delta \gamma_1=  \delta \gamma_2 =0$. Let us calculate the Laplacian of 
their wedge product:

$\triangle (\gamma_1 \wedge \gamma_2 ) = (d\delta + \delta d) (\gamma_1 \wedge \gamma_2 ) =$

$d\delta (\gamma_1 \wedge \gamma_2 ) + (\delta d) (\gamma_1 \wedge \gamma_2 )= $

$d( \delta \gamma_1 \wedge \gamma_2 \pm  \gamma_1 \wedge \delta \gamma_2)+ \delta  (d\gamma_1 
\wedge \gamma_2 \pm \gamma_1 \wedge d\gamma_2 ) =0$.

\

We have thus proved that the product of two harmonic forms is harmonic. 
It rests to prove now that if $\gamma$ is harmonic,  so it is $*\gamma$, i.e.,
if $\triangle(\gamma)=0$ then $\triangle(*\gamma)=0$.

$\triangle(*\gamma)=(\delta d + d\delta)(*\gamma)$

$= (\delta d * + d\delta *)\gamma = ((* \delta (-1)^{(n-p+1)(p-1)+1} +
d(*d(-1)^{p(n-p)+1}))\gamma$

$= (*d\delta (-1)^{(n-p+1)(p-1)+1} (-1)^{(p-1)(n-p+1)+1} + *\delta d
(-1)^{p(n-p)+1} (-1)^{(n-p)p+1}$

$= ( *(-1)^{2(n-p+1)+2}d\delta + \delta d (-1)^{2p(n-p)+2})\gamma$

$= (*(d\delta + \delta d))\gamma = *\triangle(\gamma) = 0$.

\begin{teo}
\textbf{Theorem. } Let $E$ be the basis of $ Harm^{2k}$ formed with those
elements of $H^{2k}$ that are harmonic. Then, the signature quadratic form $Q$
is block diagonal in $\{E\}$ and the signature of $M$ fulfills 

$$ \sigma(M)=
dim(Harm_+^{2k}) - dim( Harm_-^{2k})$$

where $Harm_+^{2k},  Harm_-^{2k}$ are the
eigenspaces of $*$ with $\pm1$ eigenvalues.
\end{teo}

Proof. One can combine $2k-$forms that complete one another in a  manifold of 
dimension $4k$ to get eigenvectors of  the Hodge * operator. Since, $*$
satisfies $**=1$,  $*$ has eigenvalues $\pm1$. Let eigenvectors be such 
$*\omega^+ = \omega^+$ and
$*\omega^- = -\omega^-$. Then:

\

$ \int  \omega^+ \wedge  \omega^+ =\int  \omega^+ \wedge * \omega^+ =
(\omega^+,\omega^+) > 0$

$ \int  \omega^- \wedge  \omega^- =-\int  \omega^- \wedge * \omega^-
=-(\omega^-,\omega^-) < 0$

\

Using  the fact that $\alpha \wedge *\beta = \beta \wedge *\alpha$ and that in our case the wedge 
product is commutative, we have

\

$  \int  \omega^+ \wedge  \omega^- =-\int  \omega^+ \wedge * \omega^- =-\int 
\omega^- \wedge * \omega^+=-\int  \omega^- \wedge  \omega^+= 0$

\

So,the entries of $Q$ that are outside the diagonal are all zero. 

\

These results show that the matrix of $Q$ evaluated over a basis of harmonic
forms is block diagonal, say, $D$. This  allows to alternatively define the
signature of a manifold using the following idea: for a real number $r$, we can
define its sign as $sign(r) = r/\sqrt(r^2)$. For selfadjoint diagonal matrices
we also have $sign(D) =  D(\sqrt(D^2))^{-1}$. In a general case  $sign(A) = 
A(\sqrt(A^*A))^{-1}$ which displays a diagonal filled in $\pm 1$.  Thus

$$\sigma(M) = Trace (sign(Q))$$

where $Q$ is any matrix representing the signature  quadratic form.

When some eigenvalues of a matrix are zero, they don't enter to define the signature of a given 
matrix. So, it would be nice to face up at this point the following intrigue:  Is $Q$ degenerate 
allowing the existence of null eigenvalues?

\begin{teo}
\textbf{Theorem. } $Q$ is not degenerate over $H^*$.
\end{teo}

Proof. It is enough to restrict to harmonic representatives,   which are both closed and coclosed 
with the star version also harmonic. Let us suppose now that

$\int (\gamma_1 \wedge \gamma_2) = 0, \forall \gamma_2$.

If that is true, then
we can take  $\gamma_2 = *\gamma_1$. In that case we get

$\int (\gamma_1 \wedge *\gamma_1) = (\gamma_1, \gamma_1)=0$

from which we conclude that $\gamma_1=0$, which means that $Q$ is not degenerate and that has no 
zero eigenvalues.

\

\

\begin{teo}
\textbf{Lemma. } Let $dim M= 4k$ then  the dimension of $H^{2k}$ and the signature of $M$  have the 
same parity, i.e., $b_{2k}=\sigma(M)$ mod 2.
\end{teo}

Proof: Since the harmonic forms $Harm^{2k}$ can be decomposed into the direct sum

$Harm^{2k}=Harm_+^{2k}\oplus Harm_-^{2k}$ thus $b_{2k}=dim (Harm_+^{2k})+dim (Harm_-^{2k})$

and moreover,

$\sigma(M) = dim (Harm_+^{2k})-dim (Harm_-^{2k})$

\

Therefore

$b_{2k}=dim (Harm_+^{2k})-dim (Harm_-^{2k}) + dim (Harm_-^{2k}) + dim (Harm_-^{2k}) $

$b_{2k}= \sigma(M)+2dim (Harm_-^{2k}) $

$b_{2k}- \sigma(M)=2dim (Harm_-^{2k}) $

or $b_{2k}= \sigma(M)$ mod 2.

\

Since ``to have the same parity''  is  transitive, we  can conclude

\begin{teo}
\textbf{Corollary. } Let $dim M= 4k$ then the signature of M and the Euler characteristic has the 
same parity, i.e., $\sigma(M) = \chi(M)$ mod 2.
\end{teo}

These results allow us to correct a sophism above.

\begin{teo}
\textbf{Example: clarification of a sophism on the  signature  of
$\mathbb{CP}^2$.} 
\end{teo}

We consider $\mathbb{CP}^2$ as a compact, oriented   manifold over the reals of dimension
4. It   can be decomposed as $e⁰ + e^2 + e^4$. Its Poincaré polynomial  is:

$p(t) = 1 + t^2 + t^{4}$

which means that $b_2 = b_4 = 1$. So, $H^2$ and $H^4$ both have one generator.
Thus, $Q$ is a $1 \times 1$ matrix. Since $Q$ cannot be degenerate, the only
entry of $Q$ cannot be zero. In the sophism above it was argued that because of duplication of  
terms  every $Q$ matrix
 has zeros in its diagonal and that therefore the
signature
of   $\mathbb{CP}^2$ should be  zero. This 
generalization is in general false:  when $Q$
is computed over harmonic forms, it gets a diagonal form with all eigenvalues
different than zero.

On the other hand, we have said that   $H^2$ and   $H^4$ both have one generator. This seems to be
contradictory:  $e^4$ can be parameterized in polar coordinates and so $H^4$
would be generated by $d\theta_1d\theta_2d\theta_3d\theta_4$. This would imply
that $H^2$ would have six generators instead of one: $d\theta_1d\theta_2$, 
$d\theta_1d\theta_3$,  $d\theta_1d\theta_4$, ... This trouble is resolved by
arguing that there is only one generator in $H^2$ and that the other terms
play the role  of the $dz$ of a cylinder that generates no cohomology at all. Let
us blame complex numbers for that deficit because they include rotations in
their algebraic structure. The use of harmonic forms follows:

The cell decomposition technology allows us to model $\mathbb{CP}^2$ over $e^4$
whose $H^4$ is generated by $d\theta_1d\theta_2d\theta_3d\theta_4$. We know
that $H^2$ has only one generator and that it has a harmonic representative.
Let us name the coordinates of $e^4$ in such a way  that the following form
is a representative of the only one class of harmonic 2-forms of  $\mathbb{CP}^2$:

\

$\nu = d\theta_1 d\theta_2 + d\theta_3 d\theta_4$

\

Let us check that $\nu = d\theta_1 d\theta_2 + d\theta_3 d\theta_4$ is
harmonic, i.e., that 

$\triangle \nu = ( d\delta + \delta d) (\nu) = 0 $

\

Taking $n= 4$ and
$p=2$ in  $\delta = (-1)^{n(p+1)+1} *d*  = -*d*$, we get

\

$\triangle \nu = ( d*d* + *d* d) (\nu) = ( d*d* + *d* d) (d\theta_1 d\theta_2 +
d\theta_3 d\theta_4)  $

$ \hspace{0.7cm} = ( d*d*) (d\theta_1 d\theta_2 +
d\theta_3 d\theta_4)  + (*d* d) (d\theta_1 d\theta_2 +
d\theta_3 d\theta_4) $

\

We have that $* (d\theta_1 d\theta_2 +
d\theta_3 d\theta_4) = d\theta_1 d\theta_2 +
d\theta_3 d\theta_4 $. So:

\

$\triangle \nu =  ( d*d) (d\theta_1 d\theta_2 +
d\theta_3 d\theta_4)  + (*d* d) (d\theta_1 d\theta_2 +
d\theta_3 d\theta_4) $

\

Now, $d (d\theta_1 d\theta_2 +
d\theta_3 d\theta_4) = d (1(d\theta_1 d\theta_2 +
d\theta_3 d\theta_4))= 0$. Therefore,

\

$\triangle \nu = 0$.

\

Computed over $\nu$, matrix $Q$ takes the diagonal form:

$Q = [\int \nu \wedge \nu]$
$= [\int (d\theta_1 d\theta_2 + d\theta_3 d\theta_4) \wedge (d\theta_1 d\theta_2
+ d\theta_3 d\theta_4)] = [\int 2d\theta_1 d\theta_2 d\theta_3 d\theta_4]$

$ \hspace{0.4cm} = [2c] = [c]$

This value is different than zero because the integrand is in the same class of
the volume form, so it is positive. Therefore, matrix $Q$ has one positive
eigenvalue and nothing else. As a consequence, the signature of $\mathbb{CP}^2$ 
is 1. 

\

Our result is backed by the prediction given by the Euler Characteristic $\chi(\mathbb{CP}^2) 
= \sum (-1)^kb_k$, which becomes $\chi(\mathbb{CP}^2) =  (-1)^0 + (-1)^2 + (-1)^4 = 3$. Because the 
Euler characteristic and the signature have the same parity,  the signature is predicted to  be odd, 
such as it was found.

\begin{teo}
\textbf{Example. } Let us battle with the complex hyper-torus 
$K = \mathbb{CP}^2 \times \mathbb{CP}^2$.
\end{teo}

Since the Poincaré polynomial  of $\mathbb{CP}^2 $ is $q(t) = 1 + t^2 + t^4$, 
the Poincaré polynomial of $K = \mathbb{CP}^2 \times \mathbb{CP}^2$ is:

\

$p(t) = (1 + t^2 + t^4)^2 = 1+ t^4 + t^8 + 2t^2 + 2t^4 + 2t^6 = 1+ 2t^2 +
3t^4 + 2t^6 + t^8 $

\

This polynomial predicts that $H^4$ has 3 generators, that $H^8$ has one  and
that therefore the $Q$ matrix is $3 \times 3$. The Euler characteristic reads:
$\chi(K) =  (-1)^0 + 2(-1)^2 + 3(-1)^4 + 2(-1)^6 + (-1)^8   = 9$, so the signature is expected to 
 be odd. From $\mathbb{CP}^2 $ to $K$ there is nothing special, so we can make calculations in 
whatever basis without too much trouble:
 
 The generator of $H^8$ comes from the hypertorus $e^4 \times e^4$ so,  the volume form 
in $K$ is in the class of $d\omega = 
d\theta_1d\theta_2d\theta_3d\theta_4d\phi_1d\phi_2d\phi_3d\phi_4$. On the other hand, 
 we know that $H^4$ has 3 generators, whose precedence is the following: $3t^4 = 2t^4 +   t^2 
t^2$. This means that we have two forms that come from $2e^4$ and another that comes from   $ 
e^2 \times e^2$. In consequence, let $dg_1 =d\theta_1d\theta_2d\theta_3d\theta_4$ be the first 
generator of $H^4$  and 
let the second be  $dg_2=d\phi_1d\phi_2d\phi_3d\phi_4$. And, what about the generator that comes 
from the torus   $ e^2 \times 
e^2$?  Let it be $dg_3 = d\eta_1d\eta_2d\nu_1d\nu_2$.

From this we can compute the $Q$ matrix of
the integrals of wedge products. Since we are dealing with a dimension that is multiple of 4, 
matrix $Q$ is symmetric:

 $$Q = \bordermatrix{~ & dg_1 & dg_2 & dg_3 \cr
 dg_1 & 0 & c & d  \cr
 dg_2  & c & 0 & e \cr
 dg_3  & d & e & 0 \cr}$$

\

where $d$ and $e$ are volumes of 8-dimensional shapes others than $e^4 \times e^4$ and that we were 
unable to predict at the start. By assuming that $c = d  = e$, we get:

 $$Q = \bordermatrix{~ & dg_1 & dg_2 & dg_3 \cr
 dg_1 & 0 & c & c  \cr
 dg_2  & c & 0 & c \cr
 dg_3  & c & c & 0 \cr}$$

The characteristic polynomial of this matrix is

$p(x) = -\lambda^3 +3c^2 \lambda + 2 c^3$

We verify that $-c$ is a root: $-c^3 + 3 c^3 - 2 c^3 = 0$. Hence,

$p(x) = (\lambda + c)(\lambda + c)(\lambda - 2c)$. 

This implies that the eigenvalues are $-c$, $-c$, $2c$. The signature of $K$ is, therefore,  $-1$.

\begin{teo}
\textbf{Observation. } The signature of a manifold is the number of 
positive
eigenvalues minus the number of negative eigenvalues of the matrix of the signature quadratic 
form in whatever basis. It also can  be rewritten as

$$ \sigma(M)= dim(Harm_+^{2k}) - dim(Harm_-^{2k})$$

where the subtraction refers to finite numbers.
\end{teo}

\begin{teo}
\textbf{Definition. } The \textbf{index of a Fredholm operator $P$ } is

$$ind(P)=  dim(Ker P) - dim(coKer P) =  dim(Ker P) - dim(Ker P^*)$$

where $P^*$ is the adjoint of $P$.  An operator is Fredholm when it has a kernel and cokernel of 
finite 
dimensions and its range is closed. The cokernel is the codomain quotient  the image. Example: If 
$K$ is 
compact (the image of a bounded set 
has a compact closure), $I + K$ is Fredholm (MIT, \cite{MIT00}, 
2000). 
\end{teo}

\begin{teo}
\textbf{Suspicion. } A Fredholm operator must exist  whose index is precisely the
signature of the manifold.
\end{teo}

 In the following we will construct that operator.

\section{The signature operator }

We introduce the chirality operator which together with the Dirac operator will
allow us to define the  signature operator, whose tremendous importance will be
exhibited in a suitable index theorem.

\

The following properties of the Hodge-star  $*$ operator  will be used in this
section:

\begin{enumerate}
 \item  $** = (-1)^{p(n-p)}$   when $*$ operates over  a $p$-form. Hence, over $(p+1)$-forms 
$**=(-1)^{(p+1)(n-p-1)}$  while over $(p-1)$-forms $**=(-1)^{(p-1)(n-p+1)}$.

\item Let us define the \textbf{chirality operator} $J$  over $p$-forms  by  $J=
i^{\frac{n}{2} + p(p-1)}*$. It becomes $J= i^{\frac{n}{2} + (p+1)(p)}*$ over
$(p+1)$-forms and $J= i^{\frac{n}{2} + (p-1)(p-2)}*$ over (p-1)-forms.  

\item Recall that  over $p-$forms we have

$*\delta = -**d* = (-1)^{(n-p+1)(p-1)+1}d*$

\item $*d= (-1)^{p(n-p)+1}\delta *$.
\end{enumerate}

These identities will be used to change the relative position of $*$ with respect to $d$ and 
$\delta$.

\begin{teo}
\textbf{Theorem. } If $D$ is defined as  $D=d+\delta$ then $J$ and $D$ anticommute: $JD=-DJ$, where 
 $J= i^{\frac{n}{2} + p(p-1)}*.$ \end{teo}

Proof: $JD=J(d+\delta) = Jd + J\delta $. Recall that  $d$ rises the order of a form by one unit, 
while $\delta$  lowers the order of a form by one. This must be kept in mind when operating with 
$J$: $Jd + J\delta = i^{n/2 + (p+1)(p)}*d + i^{n/2 + (p-1)(p-2)}*\delta $. Changing the relative 
position of $*$ we have:

$JD=i^{n/2 + (p+1)(p)}(-1)^{p(n-p)+1}\delta * + i^{n/2 + (p-1)(p-2)}(-1)^{(n-p+1)(p-1)+1}d*$

$=i^{n/2 + (p-1+2)(p)}(-1)^{p(n-p)+1}\delta * + i^{n/2 + (p-1)(p)-2(p-1)}(-1)^{(n-p+1)(p-1)+1}d*$.

$=i^{n/2 + (p-1)(p)}[i^{2p}(-1)^{p(n-p)+1}\delta + i^{-2(p-1)}(-1)^{(n-p+1)(p-1)+1}d]*$

$= [i^{2p}(-1)^{p(n-p)+1}\delta + i^{-2(p-1)}(-1)^{(n-p+1)(p-1)+1}d]J  $

$= [(-1)^{p+ p(n-p)+1}\delta + (-1)^{-(p-1)+(n-p+1)(p-1)+1}d]J$.

Now, let us see that if $n$ is even then the exponents in this last expression are odd. To see 
this, recall that $p(p-1)$ is always even, that if $n$ is even then so is $np$ and  that $p^2$ 
always behave like $p$, i.e., module two we have the following equalities:

$p+ p(n-p)+1 = p + np -p^2 +1 = p + 0 -p +1 =1 $(odd).

$-(p-1)+(n-p+1)(p-1)+1 = (p-1)(-1+n-p+1)+1 = (p-1)(n-p)+1 = pn -p^2 -n + p +1 = 0-p-0+p+1 = 1$(odd).

Putting all together:

$JD= - (\delta +d)J = -DJ$, and so, $J$ and $D$ anticommute.

\begin{teo}
\textbf{Example. } Let us illustrate over $\mathbb{R}^{4}$ the property $JD =
-DJ$. Variables are denoted $x,y,z,v$, an order that defines the orientation.
\end{teo}

Our definitions and identities for   $p-$forms are:

\

$\delta = -*d*$

$D = d + \delta = d -*d* $

$J= i^{\frac{n}{2} + p(p-1)}*$

\

To test the relation of $DJ $ and $JD$  we make a try with $\omega = fdxdy$,
where $f$ is a real function of the coordinates. We have:

\

 $d\omega = (\Sigma \frac{\partial f}{\partial x_i} dx_i) dxdy = \frac{\partial
f}{\partial z} dzdxdy + \frac{\partial f}{\partial v} dvdxdy $

$ = \frac{\partial f}{\partial z} dxdydz + \frac{\partial f}{\partial v} dxdydv
$

\

Since $\delta = -*d*$, we get:

\

$\delta \omega = -*d*(fdxdy) = -*d(fdzdv) = -*(\frac{\partial f}{\partial x}
dxdzdv + \frac{\partial f}{\partial y} dydzdv) $

$ \delta \omega = - \frac{\partial f}{\partial x} dy + \frac{\partial
f}{\partial y} dx $

\

To calculate $*$, we use the next trick: $(*\alpha)$ must complement $\alpha$ to
fill in  the volume form $dxdydzdv$ and the sign must be adjusted accordingly.
Examples over $\mathbb{R}^{4}$:

\

$*fdxdy = fdzdv$ because $(dxdy)(dzdv) =  dxdydzdv$.

$*fdxdzdv = fdy$ because $ dxdzdvdy = -dxdzdydv = dxdydzdv$.

$*hdydzdv = -hdx$  because

$dydzdv(-dx) = -dydzdvdx = dydzdxdv = -dydxdzdv= dxdydzdv$.

\

Plugging  results into $D = d + \delta$, we get:

\

$D \omega = d\omega  + \delta \omega = \frac{\partial f}{\partial z} dxdydz +
\frac{\partial f}{\partial v} dxdydv $
$- \frac{\partial f}{\partial x} dy + \frac{\partial f}{\partial y} dx$

\

Besides, we have that in   $\mathbb{R}^{4}$, the operator $J = i^{\frac{n}{2} +
p(p-1)}*$  takes over p-forms the next values:

\

For 0-forms: $J =  i^{\frac{4}{2} + 0}* = i^2*  = -*$

For 1-forms: $J =  i^{\frac{4}{2} + 1(1-1)}* = i^2*  = -*$

For  2-forms: $J = i^{\frac{4}{2} + 2(2-1)} = i^4  = *$

For 3-forms: $J = i^{\frac{4}{2} + 3(3-1)} = i^8 = *$

For 4-forms: $ J = i^{\frac{4}{2} + 4(4-1)} = i^{14} = -*$

\

Applying $J$ to  $D\omega$ we get:

\

$JD \omega = J (\frac{\partial f}{\partial z} dxdydz + \frac{\partial
f}{\partial v} dxdydv$
$- \frac{\partial f}{\partial x} dy + \frac{\partial f}{\partial y} dx) $

$JD \omega = *\frac{\partial f}{\partial z}dxdy dz + * \frac{\partial
f}{\partial v} dxdydv$
$- (-*) \frac{\partial f}{\partial x} dy + (-*)\frac{\partial f}{\partial y} dx
$

$JD \omega = *\frac{\partial f}{\partial z}dxdy dz + * \frac{\partial
f}{\partial v} dxdydv$
$+* \frac{\partial f}{\partial x} dy -*\frac{\partial f}{\partial y} dx $

$JD \omega = \frac{\partial f}{\partial z} dv - \frac{\partial f}{\partial v}
dz$
$- \frac{\partial f}{\partial x} dxdzdv - \frac{\partial f}{\partial y} dydzdv$

\

This result must be compared with $DJ\omega$:

\

$DJ\omega = D (J(fdxdy)  )  = D(*fdxdy ) = D (fdzdv) = (d + \delta) (fdzdv) = d
(fdzdv) +    \delta (fdzdv)  = d (fdzdv)    -*d* (fdzdv)   $

$= \frac{\partial f}{\partial x} dxdzdv + \frac{\partial f}{\partial y} dydzdv$
$  -*d (fdxdy) $

$= \frac{\partial f}{\partial x} dxdzdv + \frac{\partial f}{\partial y} dydzdv$
$  -* (\frac{\partial f}{\partial z} dzdxdy + \frac{\partial f}{\partial v}
dvdxdy) $

$= \frac{\partial f}{\partial x} dxdzdv + \frac{\partial f}{\partial y} dydzdv$
$  -* (\frac{\partial f}{\partial z} dxdydz + \frac{\partial f}{\partial v}
dxdydv) $

$= \frac{\partial f}{\partial x} dxdzdv + \frac{\partial f}{\partial y} dydzdv$
$  -* \frac{\partial f}{\partial z} dxdydz -* \frac{\partial f}{\partial v}
dxdydv $

$= \frac{\partial f}{\partial x} dxdzdv + \frac{\partial f}{\partial y} dydzdv$
$  -\frac{\partial f}{\partial z} dv + \frac{\partial f}{\partial v} dz $

\

We have verified that $JD \omega = -DJ \omega $.

\

\begin{teo}
\textbf{Theorem. } $J^2=1$ if $n=2k$
\end{teo}

Proof. We know that if $J= i^{n/2 + p(p-1)}*$  operates over $p$-forms, it produces $n-p$-forms, so 
a new application of $J$ operates over $(n-p)$-forms and hence $J$ takes the form $J= i^{n/2 + 
(n-p)(n-p-1)}*$. Thus

$J^2= i^{n/2 + (n-p)(n-p-1)}*i^{n/2 + p(p-1)}*= i^{n/2 + (n-p)(n-p-1)+n/2 + p(p-1)}**$

Recalling that over $p-$forms $**=(-1)^{p(n-p)} = (i)^{2p(n-p)}$ we get:

$J^2= i^{n/2 + (n-p)(n-p-1)+n/2 + p(p-1)}(i)^{2p(n-p)} $

$= i^{n/2 + (n-p)(n-p-1)+n/2 + p(p-1)+ 2p(n-p)} = i^{n + (n-p)(n-p-1) + p(p-1)+ 2p(n-p)} $

$=i^{n + (n-p)(n-p-1) + p(p-1)+ 2p(n-p)} =i^{n + n^2-np-n -np +p^2 +p  + p^2-p+ 2np -2p^2}$

Simplifying we obtain

$J^2 = i^{ n^2} = i^{4k^2} = (i^4)^{k^2} =1$.

\begin{teo}
\textbf{Corollary.} When $p = n-p$, $p = \frac{n}{2}$ and $J$ sends  $ \frac{n}{2}$-forms into $ 
\frac{n}{2}$-forms. If  $J$ has eigenvalues, these must be  $\pm 1$ and for the corresponding 
eigenvector $\omega$, $J\omega = \pm \omega$. In that case, if $J\omega =\lambda \omega$ then $J^2 
\omega =\lambda J\omega= \lambda^2\omega = \omega$ then $\lambda^2= 1$ and $\lambda= \pm 1$.
\end{teo}

\begin{teo}
\textbf{Lemma. } If $J\omega = \omega$ then $JD\omega = - D\omega$, and viceversa  if  $J\omega = 
-\omega$ then $JD\omega = D\omega$. So, $J$ functions as an inversion over $D\omega$ when $\omega$ 
is a positive eigenvector  of $J$, while it functions as the identity over $D\omega$ when $\omega$ 
is a negative eigenvector of J.
\end{teo}

Proof: If $J\omega = \omega$ then $DJ\omega =  D\omega$, but $DJ\omega = -JD\omega$ because these 
operators anticommute. Hence  $ -JD\omega= D\omega$ or $ JD\omega= -D\omega$. Likewise, if    
$J\omega = -\omega$ then $DJ\omega = -JD\omega = -D\omega$ or $JD\omega = D\omega$.

\begin{teo}
\textbf{Example. } Let us illustrate over $\mathbb{R}^{4}$ the fact that 
 if $J\omega = \omega$ then $JD\omega = - D\omega$, and vice-versa  if  $J\omega
= -\omega$ then $JD\omega = D\omega$. Variables are denoted $x,y,z,v$.
\end{teo}

Let us consider the 2 form $\alpha = fdxdy + f dzdv$. Recalling that for $n= 4$
and $p = 2$, $J = *$, we have:

\

$J\alpha = * fdxdy + *f dzdv = fdzdv + f dxdy$

\

So, $\alpha$ is in the positive sector of $J$. Let us verify that  $JD\alpha =
-D\alpha$.

 Borrowing results from the previous example,  we have that for $\omega =
fdxdy$:

\

\

$D \omega = d\omega  + \delta \omega = \frac{\partial f}{\partial z} dxdydz +
\frac{\partial f}{\partial v} dxdydv $
$- \frac{\partial f}{\partial x} dy + \frac{\partial f}{\partial y} dx$

$JD \omega = \frac{\partial f}{\partial z} dv - \frac{\partial f}{\partial v}
dz$
$- \frac{\partial f}{\partial x} dxdzdv - \frac{\partial f}{\partial y} dydzdv$

\

By the same token, if $\theta = fdzdv$, 
and since $D = d -\delta = d -*d*$, we get:

\

 $D\theta = (d -*d*)(\theta)= d\theta -*d*\theta $
$ =  d(fdzdv) - *d*(f dzdv)  $

$= \frac{\partial f}{\partial x} dxdzdv + \frac{\partial f}{\partial y} dydzdv $
$ - *d(fdxdy)$

$= \frac{\partial f}{\partial x} dxdzdv + \frac{\partial f}{\partial y} dydzdv $
$ - *(\frac{\partial f}{\partial z} dxdydz + \frac{\partial f}{\partial v}
dxdydv) $

$= \frac{\partial f}{\partial x} dxdzdv + \frac{\partial f}{\partial y} dydzdv $
$- \frac{\partial f}{\partial z} dv + \frac{\partial f}{\partial v} dz$

\

and

\

$JD\theta =   *\frac{\partial f}{\partial x} dxdzdv + *\frac{\partial
f}{\partial y} dydzdv $
$+*\frac{\partial f}{\partial z} dv -* \frac{\partial f}{\partial v} dz$

 $ =   \frac{\partial f}{\partial x} dy - \frac{\partial f}{\partial y} dx $
$-\frac{\partial f}{\partial z} dxdydz - \frac{\partial f}{\partial v} dxdydv$

\

Therefore

\

$D\alpha = D\omega + D\theta =  \frac{\partial f}{\partial z} dxdydz +
\frac{\partial f}{\partial v} dxdydv $
$- \frac{\partial f}{\partial x} dy + \frac{\partial f}{\partial y} dx$

+$\frac{\partial f}{\partial x} dxdzdv + \frac{\partial f}{\partial y} dydzdv $
$- \frac{\partial f}{\partial z} dv + \frac{\partial f}{\partial v} dz$

\

while

\

$JD\alpha = JD\omega + JD\theta = $
$\frac{\partial f}{\partial z} dv - \frac{\partial f}{\partial v} dz$
$- \frac{\partial f}{\partial x} dxdzdv - \frac{\partial f}{\partial y} dydzdv$

+$  \frac{\partial f}{\partial x} dy - \frac{\partial f}{\partial y} dx $
$-\frac{\partial f}{\partial z} dxdydz - \frac{\partial f}{\partial v} dxdydv$

\

We have verified that $JD\alpha = -D\alpha$, given that $J\alpha = \alpha$.

\

By contrast,  $\beta = gdxdz + gdydx$ is in the negative sector of $J$ and so 
one shall expect that $JD\beta = D\beta$.

\

We have used the signature of the manifold to infer the index  of the signature
operator because they are equal one to another. Nevertheless, it is interesting
to calculate  the index of the signature operator from raw definitions.

\begin{teo}
\textbf{Definition. } Let $E$ be the positive sector of $J$ and $F$ the negative sector. The 
\textbf{signature operator}  $D_s$ of $D$ is the restriction of $D$  to $E$ which goes onto $F$:

$$D_s = D/E: E \rightarrow F$$

\end{teo}

$D_s $ cannot be  selfadjoint because its domain is different than the codomain. Besides, we have:

$$D^*_s = D/F: F \rightarrow E$$
$$D^*_sD_s = \triangle_{+}/E: E \rightarrow E$$
$$D_sD^*_s = \triangle_{-}/F: F \rightarrow F$$

where $\triangle_{\pm}$ are the restrictions to the $\pm$ eigenspaces of $J$ of the Laplace 
operator 
$\triangle$.

\begin{teo}
\textbf{Lemma. }  The sectors of $J$ allow the decomposition of $\triangle$.
\end{teo}

In regard with the sectors of $J$, we can write:

$$\triangle = \begin{pmatrix}\triangle_+&0\cr 0&\triangle_-\end{pmatrix} =
\begin{pmatrix} D_s^*D_s 
&0 \cr 0 & D_sD_s^* \end{pmatrix} = 
$$

Because $\triangle =D^2$, this amounts to a decomposition of the Dirac operator:

$$D=\begin{pmatrix} 0&D_s^*\cr D_s&0 \end{pmatrix}$$

\

\begin{teo}
\textbf{The Signature Index Theorem}. The index of the signature operator is equal to
the signature of $M$. Thus, the signature of $M$, which is a  global topological
invariant, can be calculated analytically by means of the index of a suitable
differential operator which operates locally.
\end{teo}

\textit{Proof:} Taking $2k-$forms and the corresponding operators we have:

$Ind(D_s)=dim(KerD_s) - dim(KerD_s^*) $

Recalling that for an operator $L$ we have $ind (L)= ind(L^*)$, we can add  to the right side

$0= dim(KerD_s^*)-dim(KerD_s^*)= dim(KerD_s^*)-dim(KerD_s)$

and we have

$Ind(D_s)=  dim(KerD_s) - dim(KerD_s^*)  + dim(KerD_s^*)-dim(KerD_s)$

$=  dim(KerD_s^*)+ dim(KerD_s) - (  dim(KerD_s)+dim(KerD_s^*))$

but $dim(KerAB) = dim(KerA) + dim(KerB)$, hence

$Ind(D_s)=  dim(KerD_s^* D_s) - dim(KerD_sD_s^*)$

$=dim(Ker\triangle_+) - dim(\triangle_-)$

$=  dim (Harm_+^{2k})-dim (Harm_-^{2k})= \sigma(M)$

\begin{teo}
\textbf{Corollary. } All eigenvalues of matrix $Q$ are real when $n = 4k$.
\end{teo}

\begin{teo}
\textbf{Example. } Let us find by direct calculation the index of the signature
operator of  $S^{4}$.
Polar coordinates:  $\theta_1$, $\theta_2$,$\theta_3$,$\theta_4$ with $\rho =
1$:
\end{teo}

To begin with we must determine the positive sector of  $J = i^{\frac{n}{2} +
p(p-1)}*$. In this case, $n=4$ and $J$ takes over p-forms the next values:

\

For 0-forms: $J =  i^{\frac{4}{2} + 0}* = i^2*  = -*$

For 1-forms: $J =  i^{\frac{4}{2} + 1(1-1)}* = i^2*  = -*$

For  2-forms: $J = i^{\frac{4}{2} + 2(2-1)} = i^4  = *$

For 3-forms: $J = i^{\frac{4}{2} + 3(3-1)} = i^8 = *$

For 4-forms: $ J = i^{\frac{4}{2} + 4(4-1)} = i^{14} = -*$

\

The candidates for eigenvectors of $J$ are:

 $f + f d\theta_1d\theta_2d\theta_3d\theta_4$,
 
$g_1d\theta_1 + g_1 d\theta_2d\theta_3d\theta_4$,

$g_2d\theta_2 + g_2 d\theta_1d\theta_3d\theta_4$,

$g_3d\theta_3 + g_3 d\theta_1d\theta_2d\theta_4$,

$g_4d\theta_4 + g_4 d\theta_1d\theta_2d\theta_3$,

$h_1d\theta_1d\theta_2 + h_1 d\theta_3d\theta_4$,

$h_2d\theta_1d\theta_3 + h_2 d\theta_2d\theta_4$,

$h_3d\theta_1d\theta_4 + h_3 d\theta_2d\theta_3$.

\

Now:

\

 $J(f + f d\theta_1d\theta_2d\theta_3d\theta_4 )  = -*f -*
fd\theta_1d\theta_2d\theta_3d\theta_4  = -fd\theta_1d\theta_2d\theta_3d\theta_4
-f $,

\

$J(g_1d\theta_1 + g_1 d\theta_2d\theta_3d\theta_4) = -*g_1d\theta_1 + *g_1
d\theta_2d\theta_3d\theta_4   $

$
 = -g_1d\theta_2d\theta_3d\theta_4 - g_1 d\theta_1$,

\

$J(g_2d\theta_2 + g_2 d\theta_1d\theta_3d\theta_4) = -*g_2d\theta_2 +* g_2
d\theta_1d\theta_3d\theta_4$

$= g_2 d\theta_1d\theta_3d\theta_4 + g_2d\theta_2 $,

\

$J(g_3d\theta_3 + g_3 d\theta_1d\theta_2d\theta_4) = -*g_3d\theta_3 + *g_3
d\theta_1d\theta_2d\theta_4 $,

$ = -g_3 d\theta_1d\theta_2d\theta_4 - g_3d\theta_3 $

\

$J(g_4d\theta_4 + g_4 d\theta_1d\theta_2d\theta_3 ) = -*g_4d\theta_4 + *g_4
d\theta_1d\theta_2d\theta_3$,

$= g_4 d\theta_1d\theta_2d\theta_3  + g_4d\theta_4 $ 

\

$J(h_1d\theta_1d\theta_2 + h_1 d\theta_3d\theta_4 ) = *h_1d\theta_1d\theta_2 + *
h_1 d\theta_3d\theta_4$

$=  h_1 d\theta_3d\theta_4 + h_1d\theta_1d\theta_2 $,

\

$J(h_2d\theta_1d\theta_3 + h_2 d\theta_2d\theta_4 ) = *h_2d\theta_1d\theta_3 +
*h_2 d\theta_2d\theta_4$

$= -h_2 d\theta_2d\theta_4 - h_2d\theta_1d\theta_3 $,

\

$J(h_3d\theta_1d\theta_4 + h_3 d\theta_2d\theta_3 ) = *h_3d\theta_1d\theta_4 +
*h_3 d\theta_2d\theta_3$

$= h_3 d\theta_2d\theta_3 + h_3d\theta_1d\theta_4$.

\

We see that the positive sector of $J$ is generated by 

\

$g_2d\theta_2 + g_2 d\theta_1d\theta_3d\theta_4$

$g_4d\theta_4 + g_4 d\theta_1d\theta_2d\theta_3 $

$h_1d\theta_1d\theta_2 + h_1 d\theta_3d\theta_4 $

$h_3d\theta_1d\theta_4 + h_3 d\theta_2d\theta_3 $

\

We must apply $D = d -\delta = d - *d*$ over these forms to see whether or not
they are in the kernel of $D$:

\

$D (g_2d\theta_2 + g_2 d\theta_1d\theta_3d\theta_4) = (d - *d*)(g_2d\theta_2 +
g_2 d\theta_1d\theta_3d\theta_4 ) $

$= dg_2d\theta_2 + dg_2 d\theta_1d\theta_3d\theta_4 - *d*g_2d\theta_2 -*d*g_2
d\theta_1d\theta_3d\theta_4$

\

$= \frac{\partial g_2}{\partial \theta_1} d\theta_1d\theta_2 + \frac{\partial
g_2}{\partial \theta_3} d\theta_3d\theta_2 +\frac{\partial g_2}{\partial
\theta_4} d\theta_4d\theta_2$

$+ \frac{\partial g_2}{\partial \theta_2}d\theta_2d\theta_1d\theta_3d\theta_4 $

 $+ *dg_2d\theta_1d\theta_3d\theta_4 $

$+*dg_2 d\theta_2$ 

\

$= \frac{\partial g_2}{\partial \theta_1} d\theta_1d\theta_2 + \frac{\partial
g_2}{\partial \theta_3} d\theta_3d\theta_2 +\frac{\partial g_2}{\partial
\theta_4} d\theta_4d\theta_2$

$- \frac{\partial g_2}{\partial \theta_2}d\theta_1d\theta_2d\theta_3d\theta_4 $

 $+ *\frac{\partial g_2}{\partial \theta_2}d\theta_2d\theta_1d\theta_3d\theta_4
$

$+*(\frac{\partial g_2}{\partial \theta_1} d\theta_1d\theta_2 + \frac{\partial
g_2}{\partial \theta_3} d\theta_3d\theta_2 +\frac{\partial g_2}{\partial
\theta_4} d\theta_4d\theta_2)$

\

$= \frac{\partial g_2}{\partial \theta_1} d\theta_1d\theta_2 + \frac{\partial
g_2}{\partial \theta_3} d\theta_3d\theta_2 +\frac{\partial g_2}{\partial
\theta_4} d\theta_4d\theta_2$

$- \frac{\partial g_2}{\partial \theta_2}d\theta_1d\theta_2d\theta_1d\theta_4 $

 $- *\frac{\partial g_2}{\partial \theta_2}d\theta_1d\theta_2d\theta_3d\theta_4
$

$+*(\frac{\partial g_2}{\partial \theta_1} d\theta_1d\theta_2
- \frac{\partial g_2}{\partial \theta_3} d\theta_2d\theta_3
-\frac{\partial g_2}{\partial \theta_4} d\theta_2d\theta_4)$

\

$= \frac{\partial g_2}{\partial \theta_1} d\theta_1d\theta_2
- \frac{\partial g_2}{\partial \theta_3} d\theta_2d\theta_3
-\frac{\partial g_2}{\partial \theta_4} d\theta_2d\theta_4$

$- \frac{\partial g_2}{\partial \theta_2}d\theta_1d\theta_2d\theta_3d\theta_4 $

$-\frac{\partial g_2}{\partial \theta_2}$

 $+ \frac{\partial g_2}{\partial \theta_1} d\theta_3d\theta_4
- \frac{\partial g_2}{\partial \theta_3} d\theta_1d\theta_4 +\frac{\partial
g_2}{\partial \theta_4} d\theta_1d\theta_3$

\

This is different than zero, so $g_2d\theta_2 + g_2 d\theta_1d\theta_3d\theta_4$
is not in the kernel of $D$.

\

Let us investigate now hat happens with $g_4d\theta_4 + g_4
d\theta_1d\theta_2d\theta_3$:

\

$D(g_4d\theta_4 + g_4 d\theta_1d\theta_2d\theta_3)= (d - *d*)(g_4d\theta_4 + g_4
d\theta_1d\theta_2d\theta_3 ) $

$= dg_4d\theta_4 + dg_4 d\theta_1d\theta_2d\theta_3 - *d*g_4d\theta_4 -*d*g_4
d\theta_1d\theta_2d\theta_3 $

\

$= \frac{\partial g_4}{\partial \theta_1} d\theta_1d\theta_4 + \frac{\partial
g_4}{\partial \theta_2} d\theta_2d\theta_4 +\frac{\partial g_4}{\partial
\theta_3} d\theta_3d\theta_4$

$+ \frac{\partial g_4}{\partial \theta_4}d\theta_4d\theta_1d\theta_2d\theta_3 $

 $+ *dg_4d\theta_1d\theta_2d\theta_3 $

$-*dg_4 d\theta_4$

\

$= \frac{\partial g_4}{\partial \theta_1} d\theta_1d\theta_4 + \frac{\partial
g_4}{\partial \theta_2} d\theta_2d\theta_4 +\frac{\partial g_4}{\partial
\theta_3} d\theta_3d\theta_4$

$- \frac{\partial g_4}{\partial \theta_4}d\theta_1d\theta_2d\theta_3d\theta_4 $

 $+ *\frac{\partial g_4}{\partial \theta_4}d\theta_4d\theta_1d\theta_2d\theta_3
$

$+*(\frac{\partial g_4}{\partial \theta_1} d\theta_1d\theta_4
 + \frac{\partial g_4}{\partial \theta_2} d\theta_2d\theta_4
 +\frac{\partial g_4}{\partial \theta_3} d\theta_3d\theta_4)$

\

$= \frac{\partial g_4}{\partial \theta_1} d\theta_1d\theta_4 + \frac{\partial
g_4}{\partial \theta_2} d\theta_2d\theta_4 +\frac{\partial g_4}{\partial
\theta_3} d\theta_3d\theta_4$

$- \frac{\partial g_4}{\partial \theta_4}d\theta_1d\theta_2d\theta_3d\theta_4 $

 $- *\frac{\partial g_4}{\partial \theta_4}d\theta_1d\theta_2d\theta_3d\theta_4
$

$+\frac{\partial g_4}{\partial \theta_1} d\theta_2d\theta_3
- \frac{\partial g_4}{\partial \theta_3} d\theta_1d\theta_3
-\frac{\partial g_4}{\partial \theta_4} d\theta_1d\theta_2$

\

$= \frac{\partial g_4}{\partial \theta_1} d\theta_1d\theta_4 + \frac{\partial
g_4}{\partial \theta_2} d\theta_2d\theta_4 +\frac{\partial g_4}{\partial
\theta_3} d\theta_3d\theta_4$

$- \frac{\partial g_4}{\partial \theta_4}d\theta_1d\theta_2d\theta_3d\theta_4 $

$-\frac{\partial g_4}{\partial \theta_4}$

$+\frac{\partial g_4}{\partial \theta_1} d\theta_2d\theta_3
- \frac{\partial g_4}{\partial \theta_3} d\theta_1d\theta_3
-\frac{\partial g_4}{\partial \theta_4} d\theta_1d\theta_2$

\

Since this is different than zero, $g_4d\theta_4 + g_4
d\theta_1d\theta_2d\theta_3$ is not in the kernel of $D$.

\

On the other hand, $h_1d\theta_1d\theta_2 + h_1 d\theta_3d\theta_4 $ cannot be
in the kernel of \\ $D = d -\delta$  because $d$ rises the degree of the form by
one while $\delta$ diminishes it by one. So, $D$ would produce the sum of a
1-form and a 3-form. The same applies to
$h_3d\theta_1d\theta_4 + h_3 d\theta_2d\theta_3$.

\

We have proved that $ dim (Ker(D_s)) = 0$.

Since the index of the signature operator $D_s$ is

\

$ind(D_s) = dim (Ker(D_s)) -dim(Ker(D^*_s))$

\

to proceed further we need to calculate  $Ker(D^*_s)$.  We have:

\

$D = d -\delta = d - *d*$

\

so

\

$D^*  = d^* - \delta^*  = \delta -(*d*)^* = \delta - *^* \delta *^* $

\

Recalling that

\

 $$ *^* = (-1)^{p(n-p)} * $$

\

we get

\

$D^*  = \delta -
*\delta *$
$= \delta - **d** = \delta - d = -D$.

\

Hence $Ker(D_s)^*  = 0$ and $\triangle^* = (D^2)^* = D^2 = \triangle$ so its
index is zero. We conclude by direct calculation that

\

$ind(D_s) = dim (Ker(D_s)) -dim(Ker(D^*_s)) = 0 - 0 = 0$

\

We have verified that for $S^4$ the index of the Signature Operator  is
equal to its signature that is zero.

\begin{teo}
\textbf{Example. } Let us verify by direct calculation that the index of $D_s $
over   the 4-torus $ T^4 = S^1 \times S^1 \times S^1 \times S^1  $ is zero.
\end{teo}

Let us notice that $S^4$ and $T^4$ both have
the same leading cell $e^4$:

\

$S^4 = $ $e^o \cup e^4$

\

while

\

 $T^4 = e^0 \cup 4 e^1 \cup 6 e^2 \cup 4e^3 \cup e^4$

 \
 
Therefore, we can reuse the calculations made for $S⁴$.

\section{Acknowledgments}

For motivation, criticisms and admonitions: Sylvie Paycha, Muhammad Hassan Zarifi, Mikhail 
Malakhaltsev.

\end{document}